\documentclass[10pt]{article}
\def\vp{{\varphi}}

\newtheorem{assumption}{Assumption}[section]
  \newtheorem{definition}{Definition}[section]
  
  \newtheorem{theorem}{Theorem}[section]
  
  \newtheorem{lemma}{Lemma}[section]
  \newtheorem{proposition}{Proposition}[section]

  \newtheorem{"definition"}{"Definition"}[section]
 \newtheorem{remark}{Remark}[section]

\newcommand{\beq}{\begin{equation}}
\newcommand{\eeq}{\end{equation}}
\newcommand{\beqno}{\begin{displaymath}}
\newcommand{\eeqno}{\end{displaymath}}
\newcommand{\beqar}{\begin{eqnarray}}
\newcommand{\eeqar}{\end{eqnarray}}
\newcommand{\beqarno}{\begin{eqnarray*}}
\newcommand{\eeqarno}{\end{eqnarray*}}

\newcommand{\norm}[2]{\| {#1} \|_{{#2}}}

\newcommand{\bec}{\begin{center}}
  \newcommand{\ec}{\end{center}}

\newcommand{\abs}[1]{\vert {#1} \vert}
\newcommand{\cref}[1]{(\ref{#1})}
\newcommand{\eref}[1]{(\ref{#1})}
\newcommand{\ben}{\begin{enumerate}}
\newcommand{\en}{ \end{enumerate}}
\newcommand{\bei}{\begin{itemize}}
\newcommand{\ei}{ \end{itemize}}
\newcommand{\bed}{\begin{description}}
\newcommand{\ed}{\end{description}}
\newcommand{\bprop}{\begin{proposition}}
\newcommand{\eprop}{\end{proposition}}
\newcommand{\bdf}{\begin{definition}}
\newcommand{\edf}{\end{definition}}
\newcommand{\bth}{\begin{theorem}}
\newcommand{\eth}{\end{theorem}}
\newcommand{\blem}{\begin{lemma}}
\newcommand{\elem}{\end{lemma}}
\newcommand{\bass}{\begin{assumption}}
\newcommand{\eass}{\end{assumption}}
\newcommand{\brem}{\begin{remark}}
\newcommand{\erem}{\end{remark}}
\newcommand{\op}[1]{{\mathbf #1}}
\newcommand{\krull}[1]{\left\{ {#1} \right\}}
\newcommand{\bracket}[1]{\left[ {#1} \right]}
\newcommand{\pa}[1]{\left( {#1} \right)}
\newcommand{\map}[3]{{#1}:{#2}\rightarrow {#3}}

\newcommand{\EpXF}[3]{E^{#1}\left[\left. #2  \right| {#3} \right]}
\newcommand{\EtqX}[3]{E_{#1}^{#2}\left[ #3 \right]}
\newcommand{\EnX}[2]{E_{#1}\left[ #2 \right]}
\newcommand{\VnX}[2]{Var_{#1}\left[ #2 \right]}

\newcommand{\Ft}[1]{{\cal F}_{#1}}

\def\vp{{\varphi}}
\newcommand{\half}{\frac{1}{2}}

\newcommand{\dfdx}[2]{\frac{\partial #1}{\partial #2}}

\newcommand{\ddfdxx}[2]{\frac{{\partial}^{2} #1}{\partial {#2}^2}}
\newcommand{\ddfdxdy}[3]{\frac{{\partial}^{2} #1}{\partial {#2}\partial {#3}}}

\def\proof{\bigskip \penalty 25\noindent{\bf Proof. }}

\newcommand{\blackslug}{\hbox{\hskip 1pt
        \vrule width 4pt height 8pt depth 1.5pt\hskip 1pt}}

\def\endproof{\blackslug \bigskip}

\begin{document}

\title{Time Inconsistent Stochastic Control in Continuous Time: Theory and Examples.
\thanks{
The authors are greatly indebted to two anonymous referees, Ivar Ekeland, Ali Lazrak, Martin Schweizer, Traian Pirvu,
 Suleyman Basak, Mogens Steffensen, and Eric B\"{o}se-Wolf for very helpful comments.
}}

{
\author{ Tomas Bj{\"{o}}rk \\
\small  Department of Finance, \\
\small  Stockholm School of Economics\\
\small  tomas.bjork@hhs.se  \\
\mbox{}\\        
\and
Mariana Khapko \\
\small  Department of Management (UTSc)\\
\small   \& Rotman School of Management \\
\small  University of Toronto, \\
\small  mariana.khapko@rotman.utoronto.ca \\
\mbox{}\\     
           \and
Agatha Murgoci \\
\small  Department of Economics and Business Economics\\
\small  Aarhus University\\
\small  agatha.murgoci@econ.au.dk
}

\date{January 10, 2016}

\maketitle
\mbox{}\\
\bec
{\large A shortened version of this paper will appear   in \\
\mbox{}\\
{\em Finance and Stochastics} \\ 
\mbox{}\\ with the title \\
\mbox{}\\
``On Time Inconsistent Stochastic Control''}
\ec

\newpage
\tableofcontents
\begin{abstract}
In this paper, which is a continuation of the  discrete time paper \cite{bjo-murFSD}, we develop a theory for  continuous time stochastic control problems which, in various ways, are time inconsistent
in the sense that they do not admit a Bellman optimality principle.   We study these problems within a game theoretic framework,
and we look for Nash subgame perfect equilibrium points.
For a general controlled continuous time Markov process and a fairly general objective functional we derive an extension of the
standard Hamilton-Jacobi-Bellman  equation, in  the form of a system of non-linear equations, for
the determination for the equilibrium strategy as well as the equilibrium value function. As applications of the general theory we study non exponential discounting,  various types of mean variance problems,  a point process example, as well as a time inconsistent linear quadratic regulator. We also present a  study of time inconsistency within the framework of a general equilibrium production economy of  Cox-Ingersoll-Ross type \cite{CIR85a}.  \\
\mbox{}
\end{abstract}
\section{Introduction}\label{int}

The purpose of this paper is to study a class of stochastic control problems in continuous time,
which have the property of being time-inconsistent in the sense that they do not allow for a Bellman optimality principle.
As a consequence of this, the very concept of optimality becomes problematic, since a strategy which is optimal given a specific 
starting point in time and space, may be  non-optimal when viewed from a later date and a different state. In this paper we attack a fairly general class of time inconsistent problems by using a game-theoretic approach, so instead of searching for {\em optimal} strategies we search for subgame perfect Nash {\em equilibrium} strategies. 
The paper presents a continuous time  version of the discrete time theory developed in our previous paper \cite{bjo-murFSD}.
Since we will build heavily on the discrete time paper, the reader is referred to that paper for motivating examples and for
more detailed discussions on conceptual issues.

\subsection{Previous literature}\label{lit}
For a detailed discussion of  the game theoretic approach to time inconsistency using Nash equilibrium points
as above the reader is referred to  \cite{bjo-murFSD}. A list of some of the most important papers on the subject
is given by \cite{Bas},  \cite{Czi},   \cite{EkeLaz}, \cite{EkeMboPir}, \cite{EkePir}, \cite{Gol}, \cite{KruSmi}, \cite{mar-sol}, \cite{PelMen},  \cite{Pol}, \cite{str},  and
 \cite{VieWei09}.

All the papers above deal with particular model choices, and different authors use different methods in order to solve the problems.
To our knowledge, the present paper, which is the continuous time part of the working paper  \cite{bjo-mur09}, is the first attempt to derive a
reasonably {\em general} (albeit Markovian) theory of time inconsistent control in continuous time.
We would, however, like to stress that for the present paper we have been greatly inspired by \cite{Bas},     \cite{EkeLaz}, and \cite{EkePir}.
\subsection{Structure of the paper}\label{stru}
The structure of the paper is roughly as follows.
\bei
\item
In Section \ref{cont} we present the basic setup, and in Section \ref{cpr} 
we discuss the  concept of {\em equilibrium}.  This concept replaces, in our setting,  the optimality concept for a standard stochastic control problem, and in  Definition \ref{eqdef} we  give a precise definition of the equilibrium control and the equilibrium value function.
\item 
Since  the equilibrium concept in continuous time  is quite delicate, we build the continuous time theory on the  discrete time theory previously developed in \cite{bjo-murFSD}.  In Section \ref{chjb} we start to  study the continuous time problem by going to the limit
for a  discretized problem, and using the results from \cite{bjo-murFSD}. This leads to an extension of the standard
HJB equation to a system of equations with an embedded static optimization problem. The limiting procedure described above is done in an informal manner. It is largely heuristic, and it thus remains to clarify precisely how the derived extended HJB system is related to the precisely defined equilibrium problem under consideration. 
\item 
The needed clarification is in fact delivered in Section \ref{vt}. In Theorem \ref{vth}, which is the main theoretical result of the paper, we give a precise statement and proof of a verification theorem.  This theorem says that a solution to the extended HJB system does indeed deliver the equilibrium control and equilibrium value function to our original problem.
\item 
In Section \ref{cex} the results of Section \ref{vt} are extended to a more general reward functional.
\item
In Section \ref{tsr} we study some extensions of the theory. 
\item 
In Section \ref{eq} we prove that for every time inconsistent problem there
exists an associated standard (i.e. time consistent) control problem which in very strong sense is equivalent to the time inconsistent problem.
\item 
In Sections \ref{nedc}-\ref{lqr} we study some examples to illustrate how the theory works in a number of concrete cases.
\item
Section \ref{eqp} is devoted to a rather detailed study of a general equilibrium model for a production economy with time inconsistent preferences.
\ei

\section{The model}\label{cont}
We now turn to the formal continuous time theory. In order to present this we need some input data.
\bdf \label{ida}
The following objects are given exogenously.
\ben
\item
A drift mapping $\map{\mu}{R_+ \times R^n \times {{R^k}}}{R^n}$.
\item 
A diffusion mapping $\map{\sigma}{R_+ \times R^n \times {{R^k}}}{M(n,d)}$, where $M(n,d)$ denotes the set of all $n \times d$ matrices.
\item 
A control constraint mapping $\map{U}{R_+ \times R^n }{2^{R^k}}$
\item 
A mapping $\map{F}{R^n \times R^n}{R}$.
\item 
A mapping $\map{G}{R^n \times R^n}{R}$.
\en
\edf 
We now consider, on the time interval $[0,T]$, a controlled SDE of the form
\beq \label{cont01}
dX_t=\mu(t,X_t,u_t)dt + \sigma (t,X_t,u_t)dW_t.
\eeq
where the state process $X$ is $n$-dimensional, the Wiener process $W$ is $d$-dimensional, and the control process $u$ is $k$-dimensional, 
with the constraint $u_t \in U(t,X_t)$. 

Loosely speaking our object is to maximize, for every initial point $(t,x)$, a reward functional of the form
\beq \label{cont001}
\EnX{t,x}{F(x,X_T)}  + G\pa{x, \EnX{t,x}{X_T}}.
\eeq
This  functional is not of a form which is suitable for dynamic programming,  and it will be discussed in detail below, but first we need to specify our class of controls.
In this paper we  restrict the controls to {\bf admissible feedback control laws}.

\bdf \label{condef}
An {\bf admissible control law} is a map $\map{\op{u}}{[0,T] \times R^n}{R^k}$ satisfying the following conditions:
\ben
\item
For each $(t,x)\in [0,T] \times R^n$ we have $\op{u}(t,x) \in U(t,x)$.
\item 
For each initial point $(s,y)\in [0,T] \times R^n$ the SDE
\beq \label{cont02}
dX_t=\mu(t,X_t,\op{u}(t,X_t))dt + \sigma (t,X_t,\op{u}(t,X_t))dW_t.
\eeq
has a unique strong solution denoted by $X^{\op{u}}$.
\en
The class of admissible control laws is denoted by $\op{U}$ We will sometimes use the notation $\op{u}_t(x)$ instead of $\op{u}(t,x)$. 
\edf

We now go on to define the controlled infinitesimal generator of the SDE above.
In the present paper we use the (somewhat non-standard) convention that the infinitesimal operator acts on the time variable as well as on the space variable, so it includes the term $\dfdx{}{t}$.
The reason for this notational convention is to have formal similarity of the continuous time theory with the discrete time theory of
\cite{bjo-murFSD}. It will facilitate some arguments below considerably.

\bdf \label{cgen}
Consider the SDE (\ref{cont01}), and let $'$ denote matrix transpose.
\bei
\item
For any fixed  $u \in {{R^k}}$, the functions $\mu^u$, $\sigma^u$ and $C^u$ 
are defined by
\beqarno
\mu^u (t,x)&=&\mu (t,x,u),\\
\sigma^u (t,x)&=&\sigma (t,x,u),\\
C^u(t,x)&=&\sigma (t,x,u)\sigma (t,x,u)'.
\eeqarno
\item
For any admissible control law ${\mathbf u}$,  the functions $\mu^{\mathbf u}$,
 $\sigma^{\mathbf u}$, $C^{\mathbf u}(t,x)$  are defined by
\beqarno
\mu^{\mathbf u} (t,x)&=&\mu (t,x,{\mathbf u}(t,x)),\\
\sigma^{\mathbf u} (t,x)&=&\sigma (t,x,{\mathbf u}(t,x)),\\
C^{\mathbf u}(t,x)&=&\sigma (t,x,{\mathbf u}(t,x))\sigma (t,x,{\mathbf u}(t,x))'.
\eeqarno
\item
For any fixed  $u \in {{R^k}}$, the  operator 
${\op A}^u$ is defined by
$$
{\op A}^u=\dfdx{}{t} +\sum_{i=1}^n\mu_i^u(t,x)\dfdx{}{x_i}+
\half \sum_{i,j=1}^nC_{ij}^u(t,x)\ddfdxdy{}{x_i}{x_j}.
$$
\item
For any admissible control law ${\mathbf u}$, the   operator ${\op A}^{\mathbf u}$ is defined by
$$
{\op{A}}^{\mathbf{u}}=\dfdx{}{t} + \sum_{i=1}^n\mu_i^{\mathbf u}(t,x)\dfdx{}{x_i}+
\half \sum_{i,j=1}^nC_{ij}^{\mathbf u}(t,x)\ddfdxdy{}{x_i}{x_j}.
$$
\ei
\edf

\section{Problem formulation }\label{cpr}
In order to formulate our problem we need an objective functional. We thus consider the two  functions 
$F$ and $G$ from Definition \ref{ida}.

\bdf \label{rf}
For a fixed $(t,x) \in [0,T] \times R^n$, and a fixed  admissible control law $\op{u}$,  the  corresponding {\bf reward functional} $J$ is defined by
\beq \label{cpr1}
J(t, x,\op{u})=\EnX{t,x}{F(x,X_T^{\op{u}})}  + G\pa{x, \EnX{t,x}{X_T^{\op{u}}}}.
\eeq
\edf
\begin{remark}
In Section \ref{cex} we will consider a more general reward functional. The  restriction to 
the functional \cref{cpr1} above is done in order to minimize the notational complexity of the derivations below, 
which otherwise would be somewhat messy. 
\end{remark}
In order to have a non degenerate problem we need a formal integrability assumption.

\bass \label{bass}
We assume that for each initial point  $(t,x)\in [0,T] \times R^n$, and each admissible control law $\op{u}$, we have
$$
\EnX{t,x}{\abs{ F(x,X_T^{\op{u}})}} < \infty,\quad \EnX{t,x}{\abs{X_T^{\op{u}}}} < \infty
$$
and hence
$$
G\pa{x, \EnX{t,x}{X_T^{\op{u}}}} < \infty .
$$
\eass

Our objective is loosely that of maximizing $J(t, x,\op{u})$ for each $(t,x)$, but conceptually this turns out to be far from trivial, so instead of {\em optimal} controls we will study {\em equilibrium} controls. The equilibrium concept is made precise in Definition \ref{eqdef} below, but in order to motivate that definition we need  a brief discussion concerning the reward functional above.

We immediately note that, compared to a standard optimal control problem, the family of reward functionals above are not connected by a Bellman optimality principle. The reasons for this are as follows:
\bei
\item
The present state $x$ appears in the function $F$.
\item
In the second term we have (even apart from the appearance of the present state $x$), a nonlinear function $G$
operating on the expected value $\EnX{t,x}{X_T^{\op{u}}}$.
\ei
Since we do not have a Bellman optimality principle it is in fact unclear what we would mean by the term ``optimal'', since the optimality concept would differ at different initial times $t$ and for different initial states $x$.

The approach taken in this paper is to look at the problem from a game theoretic perspective, and look for subgame perfect Nash equilibrium points.
This will be given a precise definition below, but loosely speaking we view the game as follows:

\bei
\item
Consider a non-cooperative game, where we have one player for each point in time $t$. We refer to this player as ``Player $t$''.
\item
For each fixed $t$,  Player  $t$  can only control the process $X$  exactly at time $t$. He/she does that by choosing a control function
$\op{u}(t, \cdot)$, so the action taken at time $t$ with state $X_t$ is given by  $\op{u}(t, X_t$).
\item
Gluing together the control functions for all players we thus have a feedback control law $\map{\op{u}}{[0,T] \times R^n}{R^k}$.
\item
Given the feedback law $\op{u}$, the reward to Player   $t$ is given by the reward functional
$$
J(t, x,\op{u})=\EnX{t,x}{F(x,X_T^{\op{u}})}  + G\pa{x, \EnX{t,x}{X_T^{\op{u}}}}
$$
\ei
An informal (and slightly naive) definition of an equilibrium  for this game would be to say that  a feedback control law $\hat{\op{u}}$
is a subgame perfect Nash equilibrium if, for each $t$,  it has the following property:
\bei
\item
If for each  $s > t$, Player  $s$  chooses the control $\hat{\op{u}}(s, \cdot)$, then it is optimal  for Player  $t$ to choose
$\hat{\op{u}}(t, \cdot)$.
\ei
A definition like this works well in discrete time, but  in continuous time this is not a bona fide definition. Since Player $t$ can only choose the control $u_t$ exactly at time $t$, he only influences the control on a time set of Lebesgue measure zero, so for a controlled SDE of the form \cref{cont01} the  control chosen by an individual player will have no effect whatsoever on the dynamics of the
process. We thus need another definition of the equilibrium concept, and we will in fact follow an approach first
taken by  \cite{EkeLaz} and \cite{EkePir}. The formal definition of equilibrium is now as follows.

\bdf \label{eqdef}
Consider an admissible control law $\hat{\op{u}}$ (informally  viewed as a candidate  equilibrium law).
Choose an arbitrary admissible control law $\op{u} \in {\bf U}$ and a fixed real number $h >0$. Also fix an arbitrarily chosen initial point $(t,x)$.
Define the control law $\op{u}_h$ by
$$
\op{u}_h(s,y)=
\left \{
\begin{array}{ccl}
\op{u}(s,y),&&  \mbox{for} \ t \leq  s <t+h, \quad y \in R^n,\\
\hat{\op{u}}(s,y),&&  \mbox{for} \ t+h \leq s \leq T,\quad y \in R^n.
\end{array}
\right.
$$
If
$$
\liminf_{h \rightarrow 0} \frac{J(t, x,\hat{\op{u}})-J(t, x,\op{u}_h)}{h} \geq 0,
$$
for all $\op{u} \in {\bf U}$,
we say that $\hat{\op{u}}$ is an {\bf equilibrium} control law. Corresponding to the equilibrium law $\hat{\op{u}}$ 
we define the equilibrium {\bf value function} $V$  by
$$
V(t,x)=J(t, x,\hat{\op{u}}).
$$
\edf
We will sometimes refer to this as an {\em intrapersonal equilibrium}, since it can be viewed as a game between 
different future manifestations of you own preferences.
\begin{remark} \label{eqrem}
This is our continuous time formalization of the corresponding discrete time equilibrium concept.

Note the necessity of dividing by $h$, since for most models  we trivially would have
$$
\lim_{h \rightarrow 0} \krull{J(t, x,\hat{\op{u}})-J(t, x,\op{u}_h)}= 0.
$$

We also note that we do not get a perfect correspondence with the discrete time equilibrium concept,
since if the limit above equals zero for all $\op{u} \in {\bf U}$, it is not clear that this corresponds to a maximum or just to a stationary point.
\end{remark}
\begin{remark} 
There may exist multiple equilibria, so the equilibrium value function should strictly 
be denoted by $V(t,x, \hat{\op{u}})$ but we use $V(t,x)$ for ease of notation.
\end{remark}

\section{An informal derivation of the extended HJB equation}\label{chjb}
We now assume that there exists an equilibrium control law $\hat{\op{u}}$ (not necessarily unique) and we go
on to derive an extension of the standard Hamilton-Jacobi-Bellman (henceforth HJB) equation for the determination
of the corresponding value function $V$. To clarify the logical structure of the derivation we outline our strategy  as follows.
\bei
\item
We discretize (to some extent) the continuous time problem. We then use our results from discrete time theory to obtain a discretized
recursion for $\hat{\op{u}}$ and we then let the time step tend to zero.
\item
In the limit we obtain our continuous time extension of the HJB equation. Not surprisingly it will in fact be
a system of equations.
\item
In the discretizing and limiting procedure we mainly rely on informal heuristic reasoning. In particular we have
do {\bf not} claim that the derivation is a rigorous one. The derivation is, from a logical point of view, only of motivational value.
\item
In Section \ref{vt} we then go on to show that our (informally derived) extended HJB equation is in fact the ``correct'' one,  by proving a rigorous verification theorem.
\ei
\subsection{Deriving the equation}\label{cpd}
In this section we will, in an informal and heuristic way, derive a continuous time extension of the HJB equation.
Note again that we have no claims to rigor in the derivation, which is only  motivational. To this end we assume
that there exists an equilibrium law $\hat{\op{u}}$ and we argue as follows.
\bei
\item
Choose an arbitrary initial point $(t,x)$. Also choose a ``small'' time increment $h >0$ and an arbitrary admissible control $\op{u}$.
\item
Define the control law $\op{u}_h$ on the time interval $[t,T]$ by
$$
\op{u}_h(s,y)=
\left \{
\begin{array}{ccl}
\op{u}(s,y),&&  \mbox{for} \ t \leq  s <t+h, \quad y \in R^n,\\
\hat{\op{u}}(s,y),&&  \mbox{for} \ t+h \leq s \leq T,\quad y \in R^n.
\end{array}
\right.
$$
\item
If now $h$ is ``small enough''  we expect to have
$$
J(t,x,\op{u}_h) \leq J(t,x,\hat{\op{u}}),
$$
 and in the limit as $h \rightarrow 0$ we should have equality if $\op{u}(t,x)=\hat{\op{u}}(t,x)$.
\ei
We now refer to the  discrete time results, as well as the notation, 
from Theorem 3.13 of \cite{bjo-murFSD}, with $n$ and $n+1$ replaced by $t$ and $t+h$. We then obtain the inequality
$$
\pa{\op{A}_h^{\op{u}}V}(t,x)  -  \pa{\op{A}_h^{\op{u}}f}(t,x,x)+\pa{\op{A}_h^{\op{u}}f^x}(t,x)  -
\op{A}_h^{\op{u}} \pa{{G \diamond g}}(t,x) + \pa{\op{H}_h^{\op{u}}g}(t,x)  \leq 0
$$
Here we have used the following notation from \cite{bjo-murFSD}.
\bei
\item
For any fixed $y \in R^n$ the mapping $\map{f^y}{[0,T]\times R^n}{R}$  is defined by
$$
f^y(t,x)=\EnX{t,x}{F\pa{y,X_T^{\hat{\op{u}}}}}
$$
\item
The function $\map{f}{[0,T] \times R^n \times R^n}{R}$ is defined by
$$
f(t,x,y)=f^y(t,x).
$$
We will also, with a slight abuse of notation, denote the entire family of functions $\krull{f^y;\ y\in R^n}$ by $f$.
\item
For any function $k(t,x)$ the operator $\op{A}_h^{\op{u}}$ is defined by
$$
\pa{\op{A}_h^{\op{u}}k}(t,x)= \EnX{t,x}{k(t+h, X_{t+h}^u)}-k(t,x).
$$
\item
The function $\map{g}{[0,T] \times R^n}{R^n}$ is defined by
$$
g(t,x)=\EnX{t,x}{X_T^{\hat{\op{u}}}}.
$$
\item
The function $G \diamond g$ is defined by
$$
\pa{G \diamond g}(t,x)=G\pa{x,g(t,x)}
$$
\item 
The term $\op{H}_h^{\op{u}}g$ is defined by
$$
\pa{\op{H}_h^{\op{u}}g}(t,x)=G(x,\EnX{t,x}{g(t+h,X_{t+h}^{\op{u}})})- G(x,g(t,x)).
$$ 
\ei

We now divide the inequality by $h$ and let $h$ tend to zero.
The the operator $\op{A}_h^{\op{u}}$ will converge to the infinitesimal operator $\op{A}^u$, where $u=\op{u}(t,x)$
but the limit of  $h^{-1}\pa{\op{H}_h^{\op{u}}g}(t,x) $ requires closer investigation.

From the definition of the infinitesimal operator  we have the approximation
$$
\EnX{t,x}{g(t+h,X_{t+h}^{\op{u}})}=g(t,x)+\op{A}^{u}g(t,x) +o(h),
$$
and using a standard Taylor approximation for $G^x$ we obtain
$$
G(x,\EnX{t,x}{g(t+h,X_{t+h}^{\op{u}}})=G(x,g(t,x))+G_y(x,g(t,x))\cdot \op{A}^ug(t,x)+o(h),
$$
where
$$
G_y(x,y)=\dfdx{G}{y}(x,y).
$$
We thus obtain
$$
\lim_{h \rightarrow 0}\frac{1}{h}\pa{\op{H}_h^ug}(t,x)=G_y(x,g(t,x))\cdot \op{A}^ug(t,x).
$$
Collecting all results we arrive at our proposed extension of the HJB equation. To stress the fact that
the arguments above are largely informal we state the equation as a definition rather than as  proposition.

\bdf \label{hjb}
The extended HJB system of equations for $V$, $f$, and $g$, is defined as follows.
\ben
\item
The function $V$ is determined by 
\beqar
 \sup_{u \in U(t,x)}   \left \{ \pa{\op{A}^uV}(t,x)  -  \pa{\op{A}^uf}(t,x,x)+\pa{\op{A}^uf^x}(t,x)\right. &&  \label{hjbv}\\
- \left. \op{A}^u \pa{{G \diamond g}}(t,x) + \pa{\op{H}^ug}(t,x)\right \}  & = & 0,\quad 0 \leq t \leq T, \nonumber \\
V(T,x)&=&F(x,x)+G(x,x). \nonumber
\eeqar
\item
For every fixed $y \in R^n$ the  function  $(t,x) \longmapsto f^y(t,x)$ is defined by
\beqar
\op{A}^{\hat{\op{u}}}f^y(t,x)&=&0,\quad 0 \leq t \leq T, \label{hjbf} \\
f^y(T,x)&=&F(y,x).\nonumber 
\eeqar
\item
The function  $g$ is defined by
\beqar
\op{A}^{\hat{\op{u}}}g(t,x)&=&0, \quad 0 \leq t \leq T, \label{hjbg} \\
g(T,x)&=&x. \nonumber
\eeqar
\en
\edf
We now have some comments on the extended HJB system.
\bei
\item
The first point to notice is that we have a {\bf system} of  equations \cref{hjbv}-\cref{hjbg} 
for the simultaneous determination of
$V$, $f$ and $g$.
\item
In the expressions above, $\hat{\op{u}}$ always denotes the  control law which realizes the supremum in the first equation.
\item
The equations \cref{hjbf}-\cref{hjbg} are the Kolmogorov backward equations for the expectations
\beqarno
f^y(t,x)&=&\EnX{t,x}{F\pa{y,X_T^{\hat{\op{u}}}}},\\
g(t,x)&=&\EnX{t,x}{X_T^{\hat{\op{u}}}}.
\eeqarno
\item
In order to solve the $V$-equation  we need to know $f$ and $g$ but these are determined by the equilibrium control law
$\hat{\op{u}}$, which in turn is determined by the $sup$-part of the $V$-equation.
\item
We have used the notation
\beqarno
f(t,x,y)&=&f^y(t,x)\\
\pa{G \diamond g}(t,x)&=&G(x,g(t,x)),\\
\op{H}^ug(t,x)&=&G_y(x,g(t,x))\cdot \op{A}^ug(t,x),\\
G_y(x,y)&=&\dfdx{G}{y}(x,y).
\eeqarno
\item
The operator $\op{A}^u$ only operates on variables within parenthesis. Thus the expression $\pa{\op{A}^uf}(t,x,x)$
is interpreted as $\pa{\op{A}^uh}(t,x)$ with $h$ defined by $h(t,x)=f(t,x,x)$. In the expression $\pa{\op{A}^uf^y}(t,x)$
the operator does not act on the upper case index $y$, which is viewed as a fixed parameter. Similarly, in the expression $\pa{\op{A}^uf^x}(t,x)$,
the operator only acts on the variables $t,x$ within the parenthesis, and does not act on the upper case index $x$.
\item
In the  case when $F(x,y)$ does not depend upon $x$, and there is no $G$ term,
the problem trivializes to a standard time consistent problem. The terms $ \pa{\op{A}^uf}(t,x,x)+\pa{\op{A}^uf^x}(t,x)$
in the $V$-equation cancel, and the system reduces to the standard Bellman equation
\beqarno
\pa{\op{A}^uV}(t,x)&=&0,\\
V(T,x)&=&F(x).
\eeqarno
\item
We note that the $g$ function above appears, in a more restricted framework,  already in \cite{Bas}, \cite{EkeLaz}, and \cite{EkePir}.
\ei
\subsection{Existence and uniqueness}\label{eu}
The task of proving existence and/or uniqueness of solutions to the extended HJB system seems (at least to us)
to be technically extremely difficult. We have no idea about how to proceed so we leave it for future research.
It is thus very much an open problem.
\section{A  Verification Theorem}\label{vt}
As we have noted above, the derivation of the continuous time extension of the HJB equation in the previous section 
was very informal. Nevertheless, it seems reasonable to expect that the system in Definition \ref{hjb} will
indeed determine the equilibrium value function $V$, but so far nothing has been formally proved.
The following two conjectures are, however,  natural.
\ben
\item
Assume that there exists an equilibrium law $\hat{\op{u}}$ and that $V$ is the corresponding value function.
Assume furthermore that $V$ is in $C^{1,2}$. Define $f^y$ and $g$ by
\beqar
f^y(t,x)&=&\EnX{t,x}{F(y,X_T^{\hat{\op{u}}})}, \label{vt1}\\
g(t,x)&=&\EnX{t,x}{X_T^{\hat{\op{u}}}}.\label{vt2}
\eeqar
We then conjecture that $V$ satisfies the extended HJB system and that $\hat{\op{u}}$ realizes the supremum in the equation.
\item
Assume that $V$, $f$, and $g$ solves the extended HJB system and that the supremum in the $V$-equation is attained
for every $(t,x)$. We then conjecture that  there exists an equilibrium law $\hat{\op{u}}$, and that it is given by the maximizing $u$
in the $V$-equation. Furthermore we conjecture that  $V$ is the corresponding equilibrium value function, and $f$ and $g$ allow for the
interpretations \cref{vt1}-\cref{vt2}.
\en
In this paper we do not attempt to prove the first conjecture. Even for a
standard time consistent control problem within an SDE framework, it is well known that this is technically quite complicated,
and it typically requires the theory of viscosity solutions. We will, however, prove the second conjecture.
This obviously has the form of a verification result, and from standard theory we would expect that it can be proved with a minimum of technical complexity. We now give the precise formulation and proof of the verification theorem, but first we need to define a function space.

\bdf \label{ic}
Consider an arbitrary admissible control $\op{u} \in {\bf U}$. A function $\map{h}{R_+ \times R^n}{R}$ is said to belong to the space $L^2(X^{\op{u}})$ if it satisfies the condition
\beq \label{ic1}
\EnX{t,x}{\int_t^T\norm{h_x(s,X_s^{\op{u}})\sigma^{\op{u}}(s,X_s^{\op{u}})}{}^2ds}< \infty
\eeq
for every $(t,x)$. In this expression $h_x$ denotes the gradient of $h$ in the $x$-variable.
\edf
We can now state and prove the main result of the present paper. 
\bth[Verification Theorem] \label{vth}
Assume that (for all $y$) the functions $V(t,x)$, {$f^y(t,x)$}, $g(t,x)$, and $\hat{\op{u}}(t,x)$ have the following properties.
\ben
\item
$V$, {$f^y$}, and $g$ solves the extended HJB system in Definition \ref{hjb}.
\item 
$V(t,x)$,  and $g(t,x)$ are smooth in the sense that they are in $C^{1,2}$, and $f(t,x,y)$ is in $C^{1,2,2}$.
\item 
The function $\hat{\op{u}}$ realizes the supremum in the $V$ equation, and $\hat{\op{u}}$ is an admissible control law.
\item
$V$,  {$f^y$},  $g$, and $G \diamond g$, as well as the function $(t,x) \longmapsto f(t,x,x)$   all belong to the space $L^2(X^{\hat{\op{u}}})$.
\en
Then $\hat{\op{u}}$
is an equilibrium law, and $V$ is the corresponding equilibrium value function. Furthermore, $f$ and $g$ can be interpreted
according to \cref{vt1}-\cref{vt2}.
\eth

\proof
The proof consists of two steps:
\bei
\item
We start by showing that $f$ and $g$
have the  interpretations \cref{vt1}-\cref{vt2} and that $V$ is the value function corresponding to
 $\hat{\op{u}}$, i.e. that $V(t,x)=J(t,x,\hat{\op{u}})$.
\item
In the second step we then prove that
$\hat{\op{u}}$ is indeed an equilibrium control law.
\ei
To show that $f$ and $g$
have the  interpretations \cref{vt1}-\cref{vt2} we apply the Ito formula to the processes $f^y(s,X_s^{\hat{\op{u}}})$ and $g(s,X_s^{\hat{\op{u}}})$.
Using  \cref{hjbf}-\cref{hjbg} and the assumed integrability conditions for $f^y$ and $g$, it follows that the processes 
$f^y(s,X_s^{\hat{\op{u}}})$ and $g(s,X_s^{\hat{\op{u}}})$ are martingales, so from the boundary conditions for $f^y$ and $g$ we obtain
our desired representations of $f^y$ and $g$ as
\beqar
f^y(t,x)&=&\EnX{t,x}{F(y,X_T^{\hat{\op{u}}})}, \label{vt0011}\\
g(t,x)&=&\EnX{t,x}{X_T^{\hat{\op{u}}}}. \label{vt0012}
\eeqar
To show that $V(t,x)=J(t,x,\hat{\op{u}})$, we use the $V$ equation \cref{hjbv} to obtain:
\beqar
  \pa{\op{A}^{\hat{\op{u}}}V}(t,x)  -  \pa{\op{A}^{\hat{\op{u}}}f}(t,x,x)+\pa{\op{A}^{\hat{\op{u}}}f^x}(t,x) &&  \nonumber \\
- \op{A}^{\hat{\op{u}}} \pa{{G \diamond g}}(t,x) + \pa{\op{H}^{\hat{\op{u}}}g}(t,x)  & = & 0, \label{vt0010}
\eeqar
where
$$
\op{H}^{\hat{\op{u}}}g(t,x)=G_y(x,g(t,x))\cdot \op{A}^{\hat{\op{u}}}g(t,x).
$$
Since $f$, and $g$ satisfies \cref{hjbf}-\cref{hjbg}, we  have
\beqarno
\pa{\op{A}^{\hat{\op{u}}}f^x}(t,x)&=&0, \label{vt6}\\
\op{A}^{\hat{\op{u}}}g(t,x)&=&0,\label{vt7}
\eeqarno
so \cref{vt0010} takes the form 
\beq \label{vt10}
\pa{\op{A}^{\hat{\op{u}}}V}(t,x)  =  \pa{\op{A}^{\hat{\op{u}}}f}(t,x,x)
+\op{A}^{\hat{\op{u}}} \pa{{G \diamond g}}(t,x)
\eeq
for all $t$ and $x$.

We now apply the Ito formula to the process $V(s,X_s^{\hat{\op{u}}})$. Integrating and taking expectations gives us
$$
\EnX{t,x}{V(T,X_T^{\hat{\op{u}}})}=V(t,x)+\EnX{t,x}{\int_t^T\op{A}^{\hat{\op{u}}}V(s,X_s^{\hat{\op{u}}})ds},
$$
where the stochastic integral part has vanished because of the integrability condition $V \in L^2(X^{\hat{\op{u}}})$.
Using \cref{vt10} we thus obtain
$$
\EnX{t,x}{V(T,X_T^{\hat{\op{u}}})}=V(t,x)+
\EnX{t,x}{\int_t^T\op{A}^{\hat{\op{u}}}f(s,X_s^{\hat{\op{u}}},X_s^{\hat{\op{u}}}ds)}+
\EnX{t,x}{\int_t^T\op{A}^{\hat{\op{u}}}\pa{{G \diamond g}}(s,X_s^{\hat{\op{u}}})ds}.
$$
In the same way we obtain
\beqarno
\EnX{t,x}{\int_t^T\op{A}^{\hat{\op{u}}}f(s,X_s^{\hat{\op{u}}},X_s^{\hat{\op{u}}})ds}&=&
\EnX{t,x}{f(T,X_T^{\hat{\op{u}}},X_T^{\hat{\op{u}}})}-f(t,x,x),\\
\EnX{t,x}{\int_t^T\op{A}^{\hat{\op{u}}}\pa{{G \diamond g}}(s,X_s^{\hat{\op{u}}})ds}&=&
\EnX{t,x}{G(X_T,g(T,X_T^{\hat{\op{u}}}))}-G(x,g(t,x)).
\eeqarno
Using this and the boundary conditions for $V$, $f$, and $g$ we get
\beqarno
\EnX{t,x}{F(X_T^{\hat{\op{u}}},X_T^{\hat{\op{u}}})+G(X_T^{\hat{\op{u}}},X_T^{\hat{\op{u}}})}&=&
V(t,x)+\EnX{t,x}{F(X_T^{\hat{\op{u}}},X_T^{\hat{\op{u}}})}-f(t,x,x)\\
&+&\EnX{t,x}{G(X_T^{\hat{\op{u}}},X_T^{\hat{\op{u}}})}-G(x,g(t,x)),
\eeqarno
i.e.
\beq \label{vt15}
V(t,x)=f(t,x,x)+G(x,g(t,x)).
\eeq
Plugging \cref{vt0011}-\cref{vt0012} into \cref{vt15} we get 
$$
V(t,x)=\EnX{t,x}{F(x,X_T^{\hat{\op{u}}})}+G(x,\EnX{t,x}{X_T^{\hat{\op{u}}}})).
$$
so we obtain the desired result
$$
V(t,x)=J(t,x,\hat{\op{u}}).
$$
\vspace{5mm}

We now go on to show that $\hat{\op{u}}$ is indeed an equilibrium law, but first we need a small temporary definition.  
For any admissible control law $\op{u}$ we define $f^{\op{u}}$ and $g^{\op{u}}$ by
\beqarno
f^{\op{u}}(t,x,y)&=&\EnX{t,x}{F(y,X_T^{\op{u}})},\\
g^{\op{u}}(t,x)&=&\EnX{t,x}{X_T^{\op{u}}}.
\eeqarno
so, in particular we have $f=f^{\hat{\op{u}}}$ and $g=g^{\hat{\op{u}}}$.
For any $h >0$, and any  admissible control law $\op{u} \in {\bf U}$, we now construct the control law $\op{u}_h$ defined in Definition \ref{eqdef}.
From Lemma 3.3 and Lemma 8.8 in \cite{bjo-murFSD}, applied to the points $t$ and $t+h$ we obtain
\beqarno
J(t, x,{\op{u}_h})&=&\EnX{t,x}{J(t+h,X_{t+h}^{\op{u}_h},\op{u}_h)}\\
&-&\krull{\EnX{t,x}{ f^{\op{u}_h}(t+h,X_{t+h}^{\op{u}_h},X_{t+h}^{\op{u}_h})}- \EnX{t,x}{f^{\op{u}_h}(t+h,X_{t+h}^{\op{u}_h},x)}}\\
&-&\krull{\EnX{t,x}{G\pa{  X_{t+h}^{\op{u}_h},g^{\op{u}_h}(t+h,X_{t+h}^{\op{u}_h})}} -G\pa{x, \EnX{t,x}{g^{\op{u}_h}(t+h,X_{t+h}^{\op{u}_h})}}}.
\eeqarno
Since  $\op{u}_h=\op{u}$ on $[t,t+h]$, we have $X_{t+h}^{\op{u}_h}=X_{t+h}^{\op{u}}$, and since $\op{u}_h=\hat{\op{u}}$ on $[t+h,T]$ we have
\beqarno
J(t+h,X_{t+h}^{\op{u}_h},\op{u}_h)&=&V(t+h,X_{t+h}^{\op{u}}),\\
f^{\op{u}_h}(t+h,X_{t+h}^{\op{u}_h},X_{t+h}^{\op{u}_h})&=&f(t+h,X_{t+h}^{\op{u}},X_{t+h}^{\op{u}}),\\
f^{\op{u}_h}(t+h,X_{t+h}^{\op{u}_h},x)&=&f(t+h,X_{t+h}^{\op{u}},x),\\
g^{\op{u}_h}(t+h,X_{t+h}^{\op{u}_h})&=&g(t+h,X_{t+h}^{\op{u}}),
\eeqarno
so we obtain
\beqarno
J(t,x,\op{u}_h)&=&\EnX{t,x}{V(t+h,X_{t+h}^{\op{u}})}\\
&-&\krull{\EnX{t,x}{ f(t+h,X_{t+h}^{\op{u}},X_{t+h}^{\op{u}})}- \EnX{t,x}{f(t+h,X_{t+h}^{\op{u}},x)}}\\
&-&\krull{\EnX{t,x}{G\pa{  X_{t+h}^{\op{u}},g(t+h,X_{t+h}^{\op{u}})}} -G\pa{x, \EnX{t,x}{g(t+h,X_{t+h}^{\op{u}})}}}.
\eeqarno
Furthermore, from the $V$-equation \cref{hjbv} we have
\beqarno
  \pa{\op{A}^uV}(t,x)  -  \pa{\op{A}^uf}(t,x,x)+\pa{\op{A}^uf^x}(t,x) &&  \\
- \op{A}^u \pa{{G \diamond g}}(t,x) + \pa{\op{H}^ug}(t,x)  & \leq & 0,
\eeqarno
where we have used the notation $u=\op{u}(t,x)$.
This gives us
\beqarno
\EnX{t,x}{V(t+h,X_{t+h}^{\op{u}})}-V(t,x)  -  \krull{ \EnX{t,x}{  f(t,X_{t+h}^{\op{u}},X_{t+h}^{\op{u}})  } - f(t,x,x)  }&&\\
+\EnX{t,x}{f(t,X_{t+h}^{\op{u}},x)}-f(t,x,x) &&  \\
- \EnX{t,x}{G\pa{t+h,g(t+h,X_{t+h}^{\op{u}})}} + G(x,g(t,x))&&\\
+G\pa{x,\EnX{t,x}{g(t+h,X_{t+h}^{\op{u}})}} - G(x,g(t,x))   & \leq & o(h),
\eeqarno
or, after simplification,
\beqarno
V(t,x) &\geq& \EnX{t,x}{V(t+h,X_{t+h}^{\op{u}})}-\EnX{t,x}{f(t,X_{t+h}^{\op{u}},X_{t+h}^{\op{u}})}
+\EnX{t,x}{f(t,X_{t+h}^{\op{u}},x)}\\
&-&\EnX{t,x}{G\pa{t+h,g(t+h,X_{t+h}^{\op{u}})}} +G\pa{x,\EnX{t,x}{g(t+h,X_{t+h}^{\op{u}})}} +o(h).
\eeqarno
Combining this with the expression for $J(t,x,\op{u}_h)$ above, and the fact that (as we have proved) $V(t,x)=J(t,x,\hat{\op{u}})$, we  obtain
$$
J(t,x,\hat{\op{u}})-J(t,x,\op{u}_h) \geq o(h),
$$
so
$$
\liminf_{h \rightarrow 0} \frac{J(t,x,\hat{\op{u}})-J(t,x,\op{u}_h)}{h}\geq 0,
$$
and we are done. \endproof

\section{The general case} \label{cex}
We now turn to the most general case of the present paper, where the functional $J$ is given by
\beqar 
J(t, x,\op{u})&=&\EnX{t,x}{\int_t^T H\pa{t,x,s,X_s^{\op{u}},\op{u}_s(X_s^{\op{u}})} ds
+   F(t,x,X_T^{\op{u}})} \nonumber \\
 &+&  G\pa{t,x, \EnX{t,x}{X_T^{\op{u}}}}.\label{tdc1}
\eeqar
To study the reward functional above we  need  a slightly modified  integrability assumption.

\begin{assumption}
We assume that for each initial point  $(t,x)\in [0,T] \times R^n$, and each admissible control law $\op{u}$, we have
\beqar
\EnX{t,x}{\int_t^T \abs{H\pa{t,x,s,X_s^{\op{u}},\op{u}_s(X_s^{\op{u}})}} ds + \abs{ F(x,X_T^{\op{u}})}} &<&  \infty,\\
\EnX{t,x}{\abs{X_T^{\op{u}}}} &<& \infty.
\eeqar
\end{assumption}
The treatment of this case is very similar to the previous one, so we directly give the final result, which is the relevant extended HJB system.

\bdf \label{tdc5}
Given the objective functional \cref{tdc1} the extended HJB system  for $V$  is given by \cref{cexv1}-\cref{cexg2} below. 
\ben
\item
The function $V$ is determined by 
\beqar
 \sup_{u \in {{R^k}}}    \{ \pa{\op{A}^uV}(t,x) + H(t,x,t,x,u)
 -(\op{A}^uf)(t,x,t,x)+ (\op{A}^uf^{tx})(t,x)  &&   \label{cexv1} \\
- \left.  \op{A}^u \pa{{G \diamond g}}(t,x) + \pa{\op{H}^ug}(t,x)\right \}  &=& 0, \nonumber 
\eeqar
with boundary condition
\beq \label{cexv2}
V(T,x)=F(T,x,x)+G(T,x,x).
\eeq
\item
For each fixed $s$ and $y$, the function $f^{sy}(t,x)$ is defined by
\beqar 
\op{A}^{\hat{\op{u}}}f^{sy}(t,x)+H(s,y,t,x, \hat{\op{u}}_t(x))&=&0, \quad 0 \leq t \leq T \label{cexf1} \\
f^{sy}(T,x)&=&F(s,y,x) \label{cexf2}
\eeqar
\item
The function $g(t,x)$ is defined by
\beqar
\op{A}^{\hat{\op{u}}}g(t,x)&=&0,\quad 0 \leq t \leq T  \label{cexg1} \\
g(T,x)&=&x. \label{cexg2}
\eeqar
\en
\edf
In the definition above, $\hat{\op{u}}$ always denotes the  control law which realizes the supremum in the $V$ equation, and we have used the notation
\beqarno
f(t,x,s,y)&=&f^{sy}(t,x),\\
\pa{G \diamond g}(t,x)&=&G(t,x,g(t,x)),\\
\op{H}^ug(t,x)&=&G_y(t,x,g(t,x))\cdot \op{A}^ug(t,x),\\
G_y(t,x,y)&=&\dfdx{G}{y}(t,x,y).
\eeqarno

Also for this case we have a verification theorem. The proof is almost identical to that of Theorem \ref{vth} so we omit it.

\bth[Verification Theorem] \label{vthg}
Assume that, for all $(s,y)$, the functions $V(t,x)$, $f^{sy}(t,x)$,  $g(t,x)$,  and $\hat{\op{u}}(t,x)$ have the following properties.
\ben
\item
$V$, $f^{sy}$,  and $g$ is a solution to the extended HJB system in Definition \ref{tdc5}.
\item 
$V$, $f^{sy}$, and $g$ are smooth in the sense that they are in $C^{1,2}$.
\item 
The function $\hat{\op{u}}$ realizes the supremum in the $V$ equation, and $\hat{\op{u}}$ is an admissible control law.
\item
$V$,  $f^{sy}$,   $g$, and $G \diamond g$, as well as the function $(t,x) \longmapsto f(t,x,t, x)$ 
all belong to the space $L^2(X^{\hat{\op{u}}})$.
\en
Then $\hat{\op{u}}$
is an equilibrium law, and $V$ is the corresponding equilibrium value function. Furthermore, $f$,  and $g$ have
 the probabilistic representations
\beqar
f^{sy}(t,x)&=&\EnX{t,x}{\int_t^T H\pa{s,y,r,X_r^{\hat{\op{u}}},\hat{\op{u}}_r(X_r^{\hat{\op{u}}})}dr+F(s,y,X_T^{\hat{\op{u}}})}, \label{cexf}\\
g(t,x)&=&\EnX{t,x}{X_T^{\hat{\op{u}}}}, \quad 0 \leq t \leq T.  \label{cexg}
\eeqar
\eth
\section{Special cases and extensions}\label{tsr}
In this section we comment on possible extensions and a couple of important special cases.
\subsection{A driving point process}
In the present paper we have, for notational clearness,  confined ourselves to a pure diffusion framework. It is, however, 
very easy to extend the theory to a case where the SDE, apart from the Wiener process,  is also driven by a marked point process with Markovian characteristics.
The extended HJB system   will look exactly the same as above, but the form of the infinitesimal operator $\op{A^u}$ will of course change, and we would need to slightly modify  the integrability assumptions.
\subsection{The case when $G=0$.} \label{GO}
In the case when the term $G$ is not present, the  $V$ equation takes the form 
\beqarno
 \sup_{u \in {{R^k}}}    \{ \pa{\op{A}^uV}(t,x) + H(t,x,t,x,u)
 -(\op{A}^uf)(t,x,t,x)+ (\op{A}^uf^{tx})(t,x) \} &=& 0   \\
\eeqarno
In this case, however, it follows from the probabilistic representation of $f$ that $f(t,x,t,x)=V(t,x)$ 
so we have a cancellation in the $V$-equation. The HJB system \cref{cexv1}-\cref{cexg2} is thus replaced by the much simpler system
\beqar
 \sup_{u \in {{R^k}}}    \{ H(t,x,t,x,u) + (\op{A}^uf^{tx})(t,x) \} &=&0, \label{cexv11}\\
\op{A}^{\hat{\op{u}}}f^{sy}(t,x)+H(s,y,t,x, \hat{\op{u}}_t(x))&=&0, \label{cexf11} \\
f^{sy}(T,x)&=&F(s,y,x). \label{cexf21}
\eeqar
\subsection{Infinite horizon}
The results above can easily be extended to  the case with infinite horizon, i.e. when $T=+ \infty$. The natural reward functional will then have the form
$$
J(t, x,\op{u})= \EnX{t,x}{\int_t^{\infty}  H\pa{t,x,s,X_s^{\op{u}},\op{u}_s(X_s^{\op{u}}) }ds } 
$$
so the functions $F$ and $G$ are not present. It is easy to see that for this case we have the  extended HJB system
\beqarno
 \sup_{u \in {{R^k}}}    \{ \pa{\op{A}^uV}(t,x) + H(t,x,t,x,u)
 -(\op{A}^uf)(t,x,t,x)+ (\op{A}^uf^{tx})(t,x) \} &=&0, \\
\lim_{T \rightarrow \infty}\EnX{t,x}{V(T,X_T^{\hat{\op{u}}})}&=&0,\\
\op{A}^{\hat{\op{u}}}f^{sy}(t,x)+H(s,y,t,x, \hat{\op{u}}_t(x))&=&0,  \\
\lim_{T \rightarrow \infty}\EnX{t,x}{f^{sy}(T,X_T^{\hat{\op{u}}})}&=&0 
\eeqarno
We also have a verification theorem where the  proof is almost identical to the earlier case.
\subsection{Generalizing $H$ and $G$}
We can easily extend the result above to the case when the term $G\pa{t,x, \EnX{t,x}{X_T^{\op{u}}}}$
is replaced by
$$
G\pa{t,x, \EnX{t,x}{k\pa{X_T^{\op{u}}}}}
$$
for some function $k$.
In this case we simply define  $g$ by
$$
g(t,x)=\EnX{t,x}{k\pa{X_T^{\hat{\op{u}}}}}.
$$
The extended HJB system then looks exactly as in Definition \ref{tdc5} above, apart from the fact that the
boundary condition for $g$ is changed to
$$
g(T,x)=k(x).
$$
See \cite{kry-ste} for an interesting application.

It is  also possible to extend the $H$  term to be of the form
$$
H(t,x,s,X_s, \EnX{t,x}{b(X_s^{\op{u}})},u_s)
$$
in the integral term of the value functional. The structure of the resulting HJB system
is fairly obvious but we have omitted it since the present HJB system is, in our opinion, complicated enough as it is.
\subsection{The case with no state dependence}\label{ssc}
We see that the general extended HJB equation is quite complicated.
In many concrete cases there are, however, cancellations between different terms in the
equation. The simplest case occurs when the objective functional has the form
$$
J(t, x,\op{u})=\EnX{t,x}{F(X_T^{\op{u}})}  + G\pa{\EnX{t,x}{X_T^{\op{u}}}},
$$
so $F$ and $G$ do not depend on the present state $x$, 
and $X$ is a scalar diffusion of the form
$$
dX_t=\mu (X_t,u_t )dt + \sigma (X_t, u_t)dW_t.
$$
In this case the extended HJB equation has the form
$$
\sup_{u \in {{R^k}}}   \krull{  \op{A}^uV(t,x)  - \op{A}^u \bracket{{G \pa{g(t,x)}}}
+ G'( g(t,x)) \op{A}^ug(t,x)   }   =  0,
$$
and a simple calculation shows that
$$
-\op{A}^u \bracket{{G \pa{g(t,x)}}} + G'( g(t,x)) \op{A}^ug(t,x)= -\half \sigma^2(x,u)G''( g(t,x))g_x^2(t,x),
$$
where $g_x=\dfdx{g}{x}$.
Thus the extended HJB equation becomes
\beq \label{ssc1}
\sup_{u \in {{R^k}}}   \krull{  \op{A}^uV(t,x) -\half \sigma^2(x,u)G''( g(t,x)) g_x^2(t,x)}   =  0,
\eeq
\subsection{A scaling result} \label{sca}
In this section we derive a small  scaling result.
Let us thus consider the objective functional \cref{tdc1} above and denote, as usual, the  equilibrium control
and value function by $\hat{\op{u}}$ and $V$ respectively.
Let $\map{\vp}{R^n}{R}$ be a fixed real valued function and  consider a new  objective functional $J_{\vp}$, defined by,
\beqno
J_{\vp}(t,x,\op{u})=\vp(x)J(t,x,\op{u}),\quad  n=0,1, \ldots ,T
\eeqno
and denote the corresponding equilibrium control
and value function by $\hat{\op{u}}_{\vp}$ and $V_{\vp}$ respectively.
Since Player  $t$ is (loosely speaking) trying to maximize $J_{\vp}(t,x,\op{u})$ over $u_t$, and $\vp(x)$ is just a scaling factor which is not affected by $u_t$
the following result is intuitively obvious. The formal proof is, however, not quite trivial.

\bprop \label{sca5}
With notation  as above we have
\beqarno
V_{\vp}(t,x)&=&\vp(x)V(t,x),\\
\hat{\op{u}}_{\vp}(t,x)&=&\hat{\op{u}}(t,x).
\eeqarno
\eprop

\proof
For notational simplicity we consider the case when $J$ is of the form
\beq \label{sca1}
J(t, x,\op{u})=\EnX{t,x}{F(x,X_T^{\op{u}})}  + G\pa{x, \EnX{t,x}{X_T^{\op{u}}}}.
\eeq
The proof for the general case has exactly the same structure.

For $J$ as above we have the extended HJB system
\beqarno
 \sup_{u \in {{R^k}}}   \left \{ \pa{\op{A}^uV}(t,x)  -  \pa{\op{A}^uf}(t,x,x)+\pa{\op{A}^uf^x}(t,x)\right. &&  \\
- \left. \op{A}^u \pa{{G \diamond g}}(t,x) + G_y(x,g(t,x))\cdot \op{A}^ug(t,x)\right \}  & = & 0,\quad 0 \leq t \leq T,\\
\op{A}^{\hat{\op{u}}}f^y(t,x)&=&0, \quad 0 \leq t \leq T,\\
\op{A}^{\hat{\op{u}}}g(t,x)&=&0, \quad 0 \leq t \leq T, \\
V(T,x)&=&F(x,x)+G(x,x),\\
f^y(T,x)&=&F(y,x),\\
g(T,x)&=&x.
\eeqarno
We now recall the probabilistic interpretations
\beqarno
V(t,x)&=&\EnX{t,x}{F(x,X_T^{\hat{\op{u}}})}  + G\pa{x, \EnX{t,x}{X_T^{\hat{\op{u}}}}}\\
f(t,x,y)&=&\EnX{t,x}{F(y,X_T^{\hat{\op{u}}})}, \label{sca25}\\
g(t,x)&=&\EnX{t,x}{X_T^{\hat{\op{u}}}}.\label{sca26}
\eeqarno
and the definition
$$
\pa{G \diamond g}(t,x)=G(x,g(t,x)).
$$
From this it follows that
$$
V(t,x)=f(t,x,x) + \pa{G \diamond g}(t,x),
$$
so the first HJB equation above can be written
$$
 \sup_{u \in {{R^k}}}   \krull{ \pa{\op{A}^uf^x}(t,x)  + G_y(x,g(t,x))\cdot \op{A}^ug(t,x)}   =  0.
$$
We now turn to $J_{\vp}$, which can be written
$$
J(t, x,\op{u})=\EnX{t,x}{F_{\vp}(x,X_T^{\op{u}})}  + G_{\vp}\pa{x, \EnX{t,x}{X_T^{\op{u}}}},
$$
where
\beqarno
F_{\vp}(x,y)&=&\vp(x)F(x,y),\\
G_{\vp}(x,y)&=&\vp(x)G(x,y),
\eeqarno
and we note that
$$
\dfdx{G_{\vp}}{y}(x,y)=\vp (x)G_y(x,y)
$$
We thus obtain the HJB equation
$$
 \sup_{u \in {{R^k}}}   \krull{ \pa{\op{A}^uf_{\vp}^x}(t,x)  + \vp (x)G_y(x,g_{\vp}(t,x))\cdot \op{A}^ug_{\vp}(t,x)}   =  0,
$$
with $f_{\vp}$ and $g_{\vp}$ defined by
\beqarno
\op{A}^{\hat{\op{u}}_{\vp}}f_{\vp}^y(t,x)&=&0, \\
\op{A}^{\hat{\op{u}}_{\vp}}g_{\vp}(t,x)&=&0,  \\
f_{\vp}^y(T,x)&=&\vp(y)F(y,x),\\
g_{\vp}(T,x)&=&x.
\eeqarno
From this it follows that we can write
\beqarno
f_{\vp}(t,x,y)&=&\vp(y)f_0(t,x,y),\\
g_{\vp}(t,x)&=&g_{0}(t,x)
\eeqarno
where
\beqarno
\op{A}^{\hat{\op{u}}_{\vp}}f_{0}^y(t,x)&=&0, \\
\op{A}^{\hat{\op{u}}_{\vp}}g_{0}(t,x)&=&0,  \\
f_{0}^y(T,x)&=&F(y,x),\\
g_{0}(T,x)&=&x.
\eeqarno
and the HJB equation  has the form
$$
 \sup_{u \in {{R^k}}}   \krull{ \vp(x)\pa{\op{A}^uf_{0}^x}(t,x)  + \vp (x)G_y(x,g_{0}(t,x))\cdot \op{A}^ug_{0}(t,x)}   =  0,
$$
or, equivalently,
$$
 \sup_{u \in {{R^k}}}   \krull{ \pa{\op{A}^uf_{0}^x}(t,x)  + G_y(x,g_{0}(t,x))\cdot \op{A}^ug_{0}(t,x)}   =  0.
$$
We thus see that the system for $f$, $g$, and $\hat{\op{u}}$ is exactly the same as that for
$f_0$, $g_0$, and $\hat{\op{u}}_{\vp}$. We thus have
\beqarno
f_{\vp}(t,x,y)&=&\vp (y)f(t,x,y),\\
g_{\vp}(t,x)&=&g(t,x),\\
\hat{\op{u}}_{\vp}&=&\hat{\op{u}}.
\eeqarno
Moreover, since
$$
V_\vp (t,x)=f_{\vp}(t,x,x) + \pa{G_{\vp} \diamond g_{\vp}}(t,x),
$$
we obtain
$$
V_\vp (t,x)=\vp(x)V(t,x).
$$
\endproof

\section{An equivalent  time consistent problem}\label{eq}
The object of the present section is to provide a  link between time inconsistent and time consistent  problems. To this end we
 go back  to the general continuous time extended HJB system \cref{cexv1}-\cref{cexg2}. The $V$-equation \cref{cexv1}  reads as
\beqarno
 \sup_{u \in {\cal U}}    \{ \pa{\op{A}^uV}(t,x) + H(t,x,t,x,u)
 -(\op{A}^uh)(t,x,t,x)+ (\op{A}^uh^{tx})(t,x)  &&  \\
-  \pa{\op{A}^uf}(t,x,t,x)+\pa{\op{A}^uf^{tx}}(t,x)  
- \left.  \op{A}^u \pa{{G \diamond g}}(t,x) + \pa{\op{H}^ug}(t,x)\right \}  &=& 0, \nonumber 
\eeqarno
Let us now assume that there exists an equilibrium control law $\hat{\op{u}}$. Using $\hat{\op{u}}$ we can then construct $f$ and $g$
by solving the associated equations \cref{cexf1}-\cref{cexg2}.  We now  define  the function $K$ by
\beqarno
K(t,x,u)&=&H(t,x,t,x,u)
-(\op{A}^uh)(t,x,t,x)+ (\op{A}^uh^{tx})(t,x)  \\
&-&  \pa{\op{A}^uf}(t,x,x)+\pa{\op{A}^uf^x}(t,x)
-  \op{A}^u \pa{{G \diamond g}}(t,x) + \pa{\op{H}^ug}(t,x).
\eeqarno
 With this definition of $K$, the equation for $V$ above and its boundary condition become
\beqarno
 \sup_{u \in {\cal U}}   \krull{ \pa{\op{A}^uV}(t,x)  +K(t,x,u)}&=&0,\\
V(T,x)&=&F(x,x)+G(x,x).
\eeqarno
We now observe,  by inspection, that this is a standard HJB equation for the standard time consistent optimal control problem to maximize
\beq \label{eqfunc}
\EnX{t,x}{\int_t^TK(s,X_s,u_s)ds + F(X_T,X_T)+G(X_T,X_T)}.
\eeq
We have thus proved the following result.

\bprop \label{eqprop}
For every time inconsistent problem in the present framework there exists a standard,
time consistent optimal control problem with the following properties.
\bei
\item
The optimal value function for the standard problem coincides with the equilibrium value function for the time inconsistent problem.
\item
The optimal control for the standard problem coincides with the equilibrium control for the time inconsistent problem.
\item
The objective functional for the standard problem is given by \cref{eqfunc}.
\ei
\eprop

We immediately remark that the Proposition above is mostly of theoretical interest, and  of little ``practical'' value. The reason is of course that in order
to formulate the equivalent standard problem we need to know the equilibrium control $\hat{\op{u}}$.

Related results can be found in \cite{Bar99}, \cite{EkeLaz}, \cite{HarLai} and \cite{lutt03}. In  these papers it is proved that,
for various models where time inconsistency stems from non-exponential discounting, there exists an equivalent standard problem (with exponential discounting).

Proposition \ref{eqprop} differs from the results in the cited references above  two ways. Firstly it differs by being quite
general and not confined to a particular model.  Secondly it differs from the  results in the cited references  by having a different structure.
In other words, for the models studied in the cited papers, the equivalent problem described in  Proposition \ref{eqprop} is structurally different from
the equivalent problems presented in the cited references.  See Section 8.3 of \cite{bjo-murFSD} for a more detailed discussion of issues of this kind.
\section{Non exponential discounting} \label{nedc}
We now illustrate the theory developed above, and the first example we consider is a
fairly general case of a control problem with non exponential discounting.
\subsection{The general case} \label{genc}
Our general model is specified as follows.
\bei
\item
We consider a controlled, not necessarily time homogeneous,  Markov process in
continuous time with controlled infinitesimal generator ${\bf A}^u$.
\item
The value functional for Player  $t$ is given by
$$
J(t,c, {\bf u})=\EnX{tx}{\int_t^T\vp (t-s)H(X_s^{\bf u},u\pa{X_s^{\bf u}})ds + \vp (T-t)\Gamma \pa{X_T^{\bf u}}},
$$
where  the discounting function $\vp(t)$, the local utility function $H(x,u)$ and the final state utility function $\Gamma (x)$ are deterministic functions.
\item
We assume that the discounting function $\vp$ is non-negative and integrable over $[0, \infty )$. Without loss of generality we assume that
$$
\vp(0)=1.
$$
\ei
In the notation of Definition \ref{tdc5} we see that we have no $G$-term
so the extended HJB equation  has the form
\beqarno
 \sup_{u \in {\cal U}}   \{ \pa{\op{A}^uV}(t,x) + C(t,x,t,x,u)
-\int_t^T(\op{A}^uc^s)(t,x,t,x)ds+ \int_t^T (\op{A}^uc^{txs})(t,x)ds &&  \\
- \left. \pa{\op{A}^uf}(t,x,t,x)+\pa{\op{A}^uf^{tx}}(t,x)
\right \}  &=& 0,
\eeqarno
We recall the relations
\beqarno
f^{sy}(t,x)&=&\EnX{t,x}{F(s,y,X_T^{\hat{\op{u}}})},\quad 0 \leq t \leq T \\
f(t,x,t,x)&=&f^{tx}(t,x),\\
c^{rys}(t,x)&=&\EnX{t,x}{C\pa{r,y,s,X_s^{\hat{\op{u}}},\hat{\op{u}}(s,X_s^{\hat{\op{u}}})}}, \quad 0 \leq t \leq r.\\
c^s(t,x,t,x)&=&c^{sxt}(t,x)
\eeqarno
as well as the notational convention that the operator ${\bf A}^u$ operates only on the variables within the parenthesis,
whereas upper case indices are treated as constant parameters.
In our case the $y$ variable is not present, so we can write
\beqarno
f^{sy}(t,x)&=&f^s(t,x)\\
c^{rys}(t,x)&=&c^{rs}(t,x)
\eeqarno
and the extended HJB equation takes the form
\beqarno
 \sup_{u \in {\cal U}}   \{ \pa{\op{A}^uV}(t,x) + C(t,t,x,u)
-\int_t^T(\op{A}^uc^s)(t,t,x)ds+ \int_t^T (\op{A}^uc^{ts})(t,x)ds &&  \\
- \left. \pa{\op{A}^uf}(t,t,x)+\pa{\op{A}^uf^{t}}(t,x)
\right \}  &=& 0,
\eeqarno
where
\beqarno
f^s(t,x)&=&\EnX{t,x}{F(s,X_T^{\hat{\op{u}}})},\quad 0 \leq t \leq T,\\
f(t,t,x)&=&f^t(t,x),\\
c^{rs}(t,x)&=&\EnX{t,x}{C\pa{r,s,X_s^{\hat{\op{u}}},\hat{\op{u}}(s,X_s^{\hat{\op{u}}})}},\\
c^s(t,t,x)&=&c^{ts}(t,x).
\eeqarno
In our case we furthermore have
\beqarno
C(t,s,x,u)&=&\vp (s-t) H(x,u),\\
F(t,x)&=&\vp (T-t)\Gamma (x).
\eeqarno
so we can write
\beqarno
f^s(t,x)&=&\vp (T-s)\gamma (t,x),\\
c^{rs}(t,x)&=&\vp (s-r)h^s(t,x),\\
c^s(t,t,x)&=&\vp (s-t)h^s(t,x)
\eeqarno
where $\gamma$ and $h$ are defined by
\beqarno
\gamma (t,x)&=&\EnX{t,x}{\Gamma (X_T^{\hat{\op{u}}})},\quad 0 \leq t \leq T,\\
h^{s}(t,x)&=&\EnX{t,x}{H\pa{X_s^{\hat{\op{u}}},\hat{\op{u}}(s,X_s^{\hat{\op{u}}})}}, \quad 0 \leq t \leq s.
\eeqarno
We now easily obtain
\beqarno
C(t,t,x,u) &=&H(x,u),\\
\op{A}^uc^{rs}(t,x)&=&\vp (s-r)\op{A}^uh^{s}(t,x),\\
\op{A}^uc^{s}(t,t,x)&=&\vp (s-t)\op{A}^uh^{s}(t,x)-\vp' (t-s)h^s(t,x),\\
\op{A}^uf^{t}(t,x)&=&\vp (T-t)\op{A}^u\gamma^s(t,x),\\
\op{A}^uf(t,t,x)&=&\vp (T-t)\op{A}^u\gamma(t,x)-\vp' (T-t)\gamma (t,x).
\eeqarno
We have thus proved the following result.

\bprop \label{nedpropc}
The extended HJB system has the form
\beqarno
 \sup_{u \in {\cal U}}  \krull{ \op{A}^uV(t,x) + H(x,u)
+ \int_t^T\vp'(s-t)h^s(t,x)ds + \vp'(T-t)\gamma (t,x)
}  &=& 0,\\
V(T,x)&=&\Gamma (x),
\eeqarno
where $h^s$ and $\gamma$ are determined by
\beqarno
\op{A}^{\hat\op{u}}h^s(t,x)&=&0,\quad  0 \leq t \leq s,\\
\op{A}^{\hat\op{u}}\gamma(t,x)&=&0, \quad 0 \leq t \leq T,\\
h^{s}(s,x)&=&H\pa{x,\hat{\op{u}}(s,x)}, \\
\gamma (T,x)&=&\Gamma (x).
\eeqarno
with probabilistic interpretation
\beqarno
\gamma (t,x)&=&\EnX{t,x}{\Gamma (X_T^{\hat{\op{u}}})},\quad 0 \leq t \leq T,\\
h^{s}(t,x)&=&\EnX{t,x}{H\pa{X_s^{\hat{\op{u}}},\hat{\op{u}}(s,X_s^{\hat{\op{u}}})}}, \quad 0 \leq t \leq s.
\eeqarno
\eprop
This generalizes the corresponding results in  \cite{EkeLaz} and \cite{EkePir}, where special cases are treated in great detail.
\subsection{Optimal investment and consumption for log utility} \label{log}
In this section we study a concrete example, namely   non exponential discounting for
 a problem of investment/consumption with log utility.

We consider a market formed by a risky asset with price process $S$
and a risk free money account with price process $B$. The price dynamics are given by
\beqarno
dS_t&=& \alpha S_tdt+\sigma S_t dW_t,\\
dB_t&=&rB_tdt,
\eeqarno
where $\alpha$ and $\sigma$ are known constants,
and $r$ is the constant short rate.

Let $u_t$ be the amount of money invested in the risky asset at time $t$, and let $c_t$ be the consumption in dollars unit time.
The value $X_t$ of a self-financing
portfolio based on $S$ and $B$ will then evolve according to the  SDE
\beq \label{basdyn}
dX_t= [rX_t+\beta u_t-c_t]dt+\sigma u_tdW_t,
\eeq
where $\beta = \alpha -r$.
The controls $u$ and $c$ are of the  form $\op{u}(x)$, $\op{c}(x)$ and
our value functional is given by
$$
J(t,x,\op{u},\op{c})=\EtqX{t,x}{}{\int_t^T \vp(s-t)\ln \bracket{\op{c}(X_s^{\op{u},\op{c}})}ds+ \vp(s-t) \ln\bracket{X_T^{\op{u},\op{c}}} }.
$$
From Proposition \ref{nedpropc} we obtain the extended HJB system as
\beqarno
 \sup_{u \in {\cal U}}  \krull{ \op{A}^uV(t,x) + \ln (c)
+ \int_t^T\vp'(s-t)h^s(t,x)ds + \vp'(T-t)\gamma (t,x)
}  &=& 0,\\
V(T,x)&=& \ln (x),
\eeqarno
with probabilistic interpretation
\beqarno
\gamma (t,x)&=&\EnX{t,x}{\ln (X_T^{\hat\op{u}, \hat\op{c}})},\quad 0 \leq t \leq T,\\
h^{s}(t,x)&=&\EnX{t,x}{\ln \pa{\hat\op{c}(X_s^{\hat\op{u}, \hat\op{c}})}} \quad 0 \leq t \leq s.
\eeqarno
The first order conditions of the maximization are given by
\beqarno
\hat{c}&=&\frac{1}{V_x},\\
\hat{u}&=&-\frac{\beta}{\sigma^2}\frac{V_x}{V_{xx}}.
\eeqarno
We now make the Ansatz
$$
V(t,x)=a(t) A \ln (x)+ d(t)
$$
where $a$ and $b$ are deterministic functions of time.   We then obtain
$$
V_t=\dot{a} \ln(x)+ \dot {b}, \quad  V_x=\frac{a}{x},\quad \quad V_{xx}=-\frac{a}{x^2},
$$
so we have
$$
\hat{c}(t,x)=\frac{x}{a(t)},\quad \quad \hat{u}(t,x)=\frac{\beta}{\sigma^2}x.
$$
Given the Ansatz, the equilibrium $X$ dynamics thus take the form
$$
dX_t=X_t\bracket{r+ \frac{\beta^2}{\sigma^2} - \frac{1}{a(t)}}dt + X_t\frac{\beta}{\sigma}dW_t.
$$
and a simple calculation gives us
$$
h^s(t,x)=\EnX{t,x}{\ln \pa{\hat\op{c}(s,X_s^{\hat\op{u}, \hat\op{c}})}}=
\ln (x)+\int_t^sb(u)du - \ln (a(s)),
$$
where
$$
b(t)=r+ \half \frac{\beta^2}{\sigma^2} - \frac{1}{a(t)}
$$
and
$$
\gamma(t,x)=\EnX{t,x}{\ln \pa{X_T^{\hat\op{u}, \hat\op{c}})}}=
\ln (x)+\int_t^Tb(u)du,
$$
We can now plug these expressions into the extended HJB equation and, with $B$ denoting the primitive function of $b$,  we obtain
\beqarno
&&\dot{a} \ln (x) + \dot {d} + \ln (x) - \ln (a) + ar+ a\frac{\beta^2}{\sigma^2}-1
- a \half \frac{\beta^2}{\sigma^2}\\
&& + \ln (x) \int_t^T\vp'(s-t)ds + \int_t^T\vp'(s-t)B(s)ds-A(t)\int_t^T\vp'(s-t)ds\\
&& - \int_t^T\vp'(s-t)a(s)ds + \vp'(T-t)\ln (x) + \vp'(T-t)\bracket{B(T)-B(t)}
=0
\eeqarno
Collecting the log terms we obtain the following ODE for $a$
\beqarno
\dot{a}+\vp (T-t)+ \vp'(T-t) &=&0,\\
a(T)&=&1.
\eeqarno
This ODE can obviously be solved directly by integration. The remaining terms in the HJB equation will give us an ODE for $d$
and, given $a$,  this can also be integrated directly. The problem is thus solved, and we note that we do not have to check the PDE:s for
$h^s(t,x)$ and $\gamma (t,x)$, since  these PDE:s follow directly from the probabilistic interpretations of $h$ and $\gamma$  that we have used above.
We thus have the following result
\bprop
With notation as above, the equilibrium solution is given by
$$
V(t,x)=a(t) \ln (x) + d(t),\quad
\hat\op{c}(t,x)=\frac{x}{a(t)},\quad
\hat\op{u}(t,x)=\frac{\beta}{\sigma^2}x.
$$
where $a$ is given above.
\eprop
\subsection{Infinite horizon} \label{nedinf}
We now move to the case of infinite horizon and we restrict ourselves to the case when the controlled Markov process is time invariant.
For this case it is natural to look for a time invariant solution, i.e. to study the case when $V$ is of the form $V(t,x)=V(0,x)=V(x)$ and when
the control $\op{u}$ is of the form $\op{u}(t,x)=\op{u}(0,x)=\op{u}(x)$. We may then set $t=0$, and $\gamma =0$ in the
extended HJB system of Proposition \ref{nedpropc}. We furthermore define the function $h(t,x)$ by
$$
h(t,x)=h^{s+t}(s,x)=h^t(0,x),
$$
and, after some elementary calculations,  we have the following result.

\bprop \label{infprop}
For the time invariant case with infinite horizon, the extended HJB system has the form
\beqarno
 \sup_{u \in {\cal U}}  \krull{ \op{G}^uV(x) + H(x,u)
+ \int_0^{\infty}\vp'(t)h(t,x)dt
}  &=& 0,
\eeqarno
where $h$ is determined by
\beqarno
{\op{G}}^{\hat\op{u}}h(t,x)&=&0,\quad  t \geq 0 \\
h(0,x)&=&H\pa{x,\hat{\op{u}}(x)},
\eeqarno
with probabilistic interpretation
\beqno
h(t,x)=\EnX{0,x}{H\pa{X_t^{\hat{\op{u}}},\hat{\op{u}}(X_t^{\hat{\op{u}}})}}. \quad t \geq 0
\eeqno
The operator ${\op{G}}^{\op{u}}$ differs from ${\op{A}}^{\op{u}}$ only by the fact that the operator $\dfdx{}{t}$ in
${\op{A}}^{\op{u}}$ is replaced by $-\dfdx{}{t}$ in ${\op{G}}^{\op{u}}$.
\eprop
\begin{remark}
Since $V$ is not  a function of $t$, we obviously have $\op{G}^uV(x)=\op{A}^uV(x)$,
but for the function $h(t,x)$ the  difference between ${\op{G}}^{\op{u}}$ and ${\op{A}}^{\op{u}}$ does indeed matter.
\end{remark}

\begin{remark}
In the time consistent case of exponential discounting  $\vp(t)=e^{-\delta t}$ we have
$$
\int_0^{\infty}\vp'(t)h(t,x)dt =-\delta \EnX{0,x}{\int_0^{\infty} e^{-\delta t}H\pa{X_t^{\hat{\op{u}}},\hat{\op{u}}(X_t^{\hat{\op{u}}})} dt}
=-\delta V(x),
$$
and we have the standard HJB equation
$$
 \sup_{u \in {\cal U}}  \krull{ \op{G}^uV(x) + H(x,u)}=\delta V(x).
$$
\end{remark}

\section{Mean-variance control} \label{mv}
In this  example we will consider dynamic mean variance optimization. This is a continuous time version of a standard Markowitz investment problem, where we penalize the risk undertaken by the conditional variance. As noted in the introduction, in a Wiener driven framework this example is studied intensively in \cite{Bas}, where the authors also consider the case of multiple assets, as well as the case of a hidden Markov process (unobservable factors) driving the
parameters of the asset price dynamics.  For illustrative purposes we first consider the simplest possible case of a Wiener driven single risky asset and re-derive the corresponding results of \cite{Bas}. We then extend the model in \cite{Bas} and study the case when the risky asset is driven by a point process as well as by a Wiener process. It should be noted that both models under consideration are special cases of the very general mean-variance theory developed in \cite{Czi}. We thus make no claim of originality for these models. They are merely intended as illustrative examples.
\subsection{The simplest case} \label{mvs}
We consider a market formed by a risky asset with price process $S$
and a risk free money account with price process $B$. The price dynamics are given by
\beqarno
dS_t&=& \alpha S_tdt+\sigma S_t dW_t,\\
dB_t&=&rB_tdt,
\eeqarno
where $\alpha$ and $\sigma$ are known constants,
and $r$ is the constant short rate.

Let $u_t$ be the amount of money invested in the risky asset at time $t$. The value $X_t$ of a self-financing
portfolio based on $S$ and $B$ will then evolve according to the  SDE
\beq 
dX_t= [rX_t+(\alpha-r)u_t]dt+\sigma u_tdW_t.
\eeq
Our value functional is given by
$$
J(t,x,\op{u})=\EtqX{t,x}{}{X_T^\op{u}}-\frac{\gamma}{2}Var_{t,x} \pa{X_T^\op{u}},
$$
so we want to maximize expected return with a penalty term for risk.
Remembering the definition for the conditional variance
$$
Var_{t,x} [X_T]=E_{t,x}[X^{2}_T]-E_{t,x}^2[X_T],
$$
we can re-write our objective functional as
$$
J(t,x,\op{u})=\EtqX{t,x}{}{F(X_T^\op{u})}-G(\EtqX{t,x}{}{X_T^\op{u}})
$$
where $F(x)= x -\frac{\gamma}{2}{x^2}$ and $G(x)=\frac{\gamma}{2}{x^2}$.
As seen in the previous sections, the term $G(\EtqX{t,x}{}{X_T})$
leads to a time inconsistent game theoretic problem.

The extended HJB equation is then given by the following PDE system:
\beqarno
\sup_u\krull{[rx+(\alpha-r)u]V_x+\half\sigma^{2}u^{2}V_{xx}-
\mathcal{A}^{u}(G\circ g)+\op{H}^{u}g}&=&0,\\
V(T,x)&=&x,\\
\mathcal{A}^{\hat{u}}g&=&0,\\
g(T,x)&=& x,
\eeqarno
where lower case index denotes the  corresponding partial derivative. This case is covered in Section \ref{ssc}, and from \cref{ssc1}
we can simplify to
\beqarno
V_t +\sup_u\krull{[rx+(\alpha-r)u]V_x+\half\sigma^{2}u^{2}V_{xx}-\frac{\gamma}{2}\sigma^{2}u^{2}g^{2}_{x}}&=&0\\
V(T,x)&=&x\\
\mathcal{A}^{\hat{u}}g&=&0\\
g(T,x)&=& x
\eeqarno

Given the linear structure of the dynamics, as well as of the boundary conditions,  it is natural to
 make the {\em Ansatz}
\beqarno
V(t,x)&=&A(t)x+B(t)\\
g(t,x)&=& a(t)x+b(t).
\eeqarno

With this trial solution the HJB equation becomes
\beqar
A_t x + B_t +\sup_u\krull{[rx+(\alpha-r)u]A-\frac{\gamma}{2}\sigma^{2}u^{2}a^{2}}&=&0,\label{HJBMV1}\\
a_t x+ b_t + [rx+(\alpha-r)\hat{u}] a &=& 0,\label{HJBMV2}\\
A(T)&=&1, \nonumber \\
B(T)&=&0, \nonumber \\
a(T)&=&1,\nonumber \\
b(T)&=&0.\nonumber
\eeqar
We first solve  the static problem embedded in \cref{HJBMV1}.
From the first order condition, we obtain the optimal control as
$$
\hat{u}(t,x) = \frac{1}{\gamma}\frac{\alpha-r}{\sigma^2}\frac{A(t)}{a^2(t)},
$$
so the optimal control does not depend on $x$.
Substituting this expression for $\hat{u}$ into   \eref{HJBMV1} we obtain:
$$
A_t x + B_t +Ar x  + \frac{1} {2\gamma} \frac {( \alpha - r )^2} {\sigma^2} \frac{ A^2 } {a^2}=0.
$$
By separation of variables we then get the following system of ODE:s.
\beqarno
A_t  +Ar &=&0,\label{ODEMV1}\\
A(T)&=&1,\nonumber\\
B_t   + \frac{1} {2\gamma} \frac {( \alpha - r )^2} {\sigma^2} \frac{ A^2 } {a^2} &=&0,\\
B(T)&=&0.\nonumber
\eeqarno
We immediately obtain
$$
A(t)= e^{r(T-t)}.
$$
Inserting this expression for $A$ into the second ODE yields
\beqar
B_t   + \frac{1} {2\gamma} \frac {( \alpha - r )^2} {\sigma^2} \frac{ e^{2r(T-t)} } {a^2}&=&0,\label{ODEB}\\
B(T)&=&0.\nonumber
\eeqar
This equation contain the unknown function $a$, and to determine this we use equation
\eref{HJBMV2}. Inserting the expression for $\hat{u}$ into \eref{HJBMV2} gives us
\beqarno
a_t x+ b_t + rx a + \frac{1}{\gamma} \frac{(\alpha - r)^2} {\sigma^2} \frac{ e^{r(T-t)} } {a} &=& 0,\\
a(T)&=&1,\\
b(T)&=&0.
\eeqarno
Again we have separation of variables and obtain the system
\beqarno
a_t  + a r &=&0,\label{ODEMV2}\\
b_t  + \frac{1}{\gamma} \frac{(\alpha - r)^2} {\sigma^2} \frac{ e^{r(T-t)} } {a} &=& 0.
\eeqarno
This yields
$$
a(t)= e^{r(T-t)},
$$
and the ODE for $b$ then takes the form
\beqarno
b_t&=&\frac{1}{\gamma} \frac{(\alpha - r)^2} {\sigma^2},\\
 b(T)&=&0.
\eeqarno
We thus have
$$
b(t)= \frac{1}{\gamma} \frac{(\alpha - r)^2} {\sigma^2}(T-t).
$$
Introducing the results in the optimal control formula, we get
$$
\hat{u}(t,x) =\frac{1}{\gamma}\frac{\alpha-r}{\sigma^2}e^{-r(T-t)}.
$$
Using the expression for $a$ above, we can go back to  equation \eref{ODEB} which now takes the form
$$
B_t  + \frac{1} {2\gamma} \frac {( \alpha - r )^2} {\sigma^2} =0,
$$
so
$$
B(t)= \frac{1} {2\gamma} \frac {( \alpha - r )^2} {\sigma^2}(T-t).
$$
Thus, the optimal value function is given by
$$
V(t,x) = e^{r(T-t)}x + \frac{1} {2\gamma} \frac {( \alpha - r )^2} {\sigma^2}(T-t)
$$
We summarize the results as follows:
\bprop \label{mvprop}
For the model above we have the following results.
\bei
\item
The optimal amount of money invested in a stock is given by
$$
\hat{\op{u}}(t,x) =\frac{1}{\gamma}\frac{\alpha-r}{\sigma^2}e^{-r(T-t)}.
$$
\item
The equilibrium value function is given by
$$
V(t,x) = e^{r(T-t)}x + \frac{1} {2\gamma} \frac {( \alpha - r )^2} {\sigma^2}(T-t).
$$
\item
The expected value of the optimal portfolio is given by
$$
\EnX{t,x}{X_T} = e^{r(T-t)}x + \frac{1}{\gamma} \frac{(\alpha - r)^2} {\sigma^2}(T-t).
$$

\ei
\eprop

Using Proposition \ref{eqprop} we can also construct the
equivalent standard time consistent optimization problem.
An easy calculation gives us the following result.

\bprop
The equivalent (in the sense of Proposition \ref{eqprop}) time consistent problem is
to maximize the functional
$$
\max_u \EtqX{t,x}{}{X_T-\frac{\gamma \sigma^2}{2}\int_t^Te^{2r(T-s)}u_s^2ds}
$$
given the dynamics \cref{basdyn}.
\eprop
We note in passing that
$$
\sigma^2u_t^2dt=d\langle X \rangle_t.
$$

\subsection{A point process extension}
We will now present an extension of the mean variance problem when the stock dynamics are driven by a jump diffusion.
We consider a single risky asset with price  $S$ and a bank account with price process $B$. The results below can be easily extended to the case of multiple assets, but for ease of exposition, we restrict ourselves to the scalar case.
The dynamics are given by
\beqarno
dS_t&=&\alpha(t,S_t) S_t dt +\sigma(t,S_t) S_t dW_t + S_{t-}\int_{\cal Z}\beta(z) \mu(dz,dt), \\
dB_t &=& rdt.
\eeqarno
Here $W$ is a scalar Wiener process and $\mu$ is a marked point process on the  mark space $\cal Z$ with
deterministic intensity measure $\lambda(dz)$.
Furthermore,  $\alpha(t,s)$, $\sigma(t,s)$ and $\beta(z)$ are known deterministic functions and  $r$ is a known constant.

As before $u_t$ denotes the amount of money invested in the stock at time $t$, and $X$
is the value process for a self financing portfolio based on $S$ and $B$.
The dynamics of $X_t$ are then given by
$$
dX_t= [rX_t+(\alpha(t,S_t,Y_t)-r)u]dt+\sigma(t, S_t,Y_t) u dW_t+ u_{t-}\int_{Z}\beta(z) \mu(dz,dt).
$$
Again we study the case of mean-variance utility, i.e.
$$
J(t,x, \op{u})= \EtqX{t,x}{}{X_T^\op{u}}-\frac{\gamma}{2}Var_{t,x} \pa{X_T^\op{u}}.
$$
The extended HJB system now has the form
\beqar
\sup_u\krull{\mathcal{A}^{u}V(t,x,s)- \mathcal{A}^{u}(G\circ g)(t,x,s)+\pa{\op{H}^{u}g}(t,x,s)}&=&0,\label{HJBmpp1}\\
V(T,x,s)&=&x,\nonumber \\
\mathcal{A}^{\hat{u}}g&=&0,\label{HJBmpp2}\\
g(T,x,s)&=& x.\nonumber
\eeqar
As before, we make the {\em Ansatz}
\beqarno
V(t,x,s)&=&A(t)x+B(t,s),\\
g(t,x,s)&=&a(t)x + b(t,s),\\
A(T) &=& 1,\\
B(T,s) &=& 0\\
a(T) &=& 1,\\
b(T,s) &=& 0.
\eeqarno
After some simple but tedious  calculations, equation \cref{HJBmpp1} can be re-written as
\beqar
&&\sup_{u}\{A_tx + B_t + A\bracket{rx+ (\alpha-r)u}+\alpha s B_s+\half\sigma^{2}s^{2}B_{ss}+A u \int_{Z}\beta(z)\lambda(dz) \nonumber \\
&&+ \int_{Z}[\underbrace{B(t,s(1+\beta(z)))- B(t,s)}_{\Delta_{\beta}B(t,s,z)}]\lambda(dz)  -
\half \sigma^2 \gamma [ au+ b_{s}s]^2  \nonumber \\
&&- \frac{\gamma}{2}\int_{Z}[au\beta(z)+\underbrace{b(t,s(1+\beta(z)))- b(t,s)}_{\Delta_{\beta}b(t,s,z)}]^2\lambda(dz)\}= 0 \label{HJBmpp201}
\eeqar

First, we solve the embedded static problem in \cref{HJBmpp201}
\beqno
\max_u \krull{ (\alpha-r)Au+ A u\int_{Z}\beta(z) \lambda(dz)-
\half \sigma^2[ au+ b_{s}s]^2-\frac{\gamma}{2}\int_{Z}[au\beta(z)+\Delta_{\beta}b] \lambda(dz)}
\eeqno
and obtain the optimal control \\
\mbox{}\\
$$
\hat{\bf u}(t,x,s) = \frac{\bracket{\alpha (t,s)-r+ \int_{Z}\beta(z) \lambda(dz)} A(t)}
{\gamma a^{2}(t)[\sigma^{2}(t,s)+\int_{Z}\beta^2(z) \lambda(dz)]}-
\frac{\sigma (t,s) b_{s}(t,s)s+ \int_{Z}\beta(z) \Delta_{\beta}b(t,s,z)\lambda(dz)}{ a(t)[\sigma^{2}(t,s)+\int_{Z}\beta^2(z)
\lambda(dz)]}
$$
\\
\mbox{}\\
Again we see that the optimal control does not depend on $x$.
We can plug the optimal control into  equation \cref{HJBmpp201} and
as before, we can separate variables to obtain an ODE for $A(t)$ and a PIDE for $B(t,s)$.
The ODE for $A$ is
\beqarno
A_t + r A = 0\\
A(T) = 0
\eeqarno
with solution $A(t)= e^{r(T-t)}$.
The PIDE for $B(t,s)$ becomes
\beqar
&&B_t +  (\alpha-r)\hat{u}+\alpha s B_s+\half\sigma^{2}s^{2}B_{ss}+A \hat{u} \int_{Z}\beta(z)\lambda(dz)
\\ &&+ \int_{Z}[\Delta_{\beta}B]\lambda(dz)  - \half \sigma^2 \gamma [ a(t)\hat{u}+ b_{s}s]^2
\\&&- \frac{\gamma}{2}\int_{Z}[a(t)\hat{u}\beta(z)+\Delta_{\beta}b]^2\lambda(dz)=  0 \label{Bpide1}\\
&&B(T,s) = 0  \label{Bpide2}
\eeqar
In order to solve this we need to determine the functions $a(t)$ and $b(t,s)$. To this end
we use \cref{HJBmpp2}.  This can be rewritten as
\beqar
&&a_t x + b_t+[rx + (\alpha -r)\hat{u}]a+\alpha s b_s+\half\sigma^{2}s^{2}b_{ss} \nonumber\\
&&+a \hat{u}\int_{Z}\beta(z)\lambda(dz)+\int_{Z}\Delta_{\beta}b\lambda(dz)= 0. \label{HJBmpp202}
\eeqar
with the appropriate boundary conditions for $a$ and $b$.
By separation of variables we obtain the  ODE
\beqarno
a_t + r a = 0,\\
a(T) = 1,
\eeqarno
and the PIDE
\beqarno
b_t+(\alpha -r)\hat{u}a+\alpha s b_s+\half\sigma^{2}s^{2}b_{ss}
+a \hat{u}\int_{Z}\beta(z)\lambda(dz)+\int_{Z}\Delta_{\beta}b\lambda(dz)&=& 0, \\
b(T,s)&=&0.
\eeqarno

From the ODE we have $a(t)= e^{r(T-t)}$ and, plugging this expression into the previous formula for $\hat{u}$, gives us
$$
\hat{u}=\frac{\alpha-r+ \int_{Z}\beta(z) \lambda(dz)}{\gamma [\sigma^{2}+\int_{Z}\beta^2(z) \lambda(dz)]}e^{-r(T-t)}
-\frac{\sigma b_{s}s+ \int_{Z}\beta(z) \Delta_{\beta}b\lambda(dz)}{ [\sigma^{2}+\int_{Z}\beta^2(z) \lambda(dz)]}e^{-r(T-t)}
$$

We can now insert this expression, as well as the formula for $a$, into the PIDE for $b$ above to obtain the PIDE

\beqarno
&&b_t+\left[\alpha-\frac{\bracket{\alpha-r+\int_{Z}\beta(y) \lambda(dy)}\sigma^2}{[\sigma^{2}+\int_{Z}\beta^2(y)
\lambda(dy)]}\right] s b_s + \half \sigma^2 s^2 b_{ss}+\frac{\bracket{\alpha-r+ \int_{Z}\beta(y) \lambda(dy)}^2 }
{\gamma [\sigma^{2}+\int_{Z}\beta^2(y) \lambda(dy)]}\\
&&\int_Z \Delta_{\beta}b(t,s,z)\krull{ 1-\frac{\alpha-r+\int_{Z}\beta(y) \lambda(dy)}{[\sigma^{2}+\int_{Z}\beta^2(y) \lambda(dy)]}
\beta(z) }\lambda(dz)= 0\\
&&b(T,s)=0
\eeqarno
This rather forbidding looking equation cannot in general be solved explicitly, but by applying a Feynman-Kac
representation theorem we can represent the solution as
\beq
b(t,s)= \EtqX{t,s}{Q}{\int_{t}^{T}\frac{\bracket{\alpha (\tau,S_{\tau})-r+ \int_{Z}\beta(z) \lambda(dz)}^2 }
{\gamma [\sigma^{2}(\tau,S_{\tau})+\int_{Z}\beta^2(z) \lambda(dz)]}d\tau}. \label{b_MPP2}
\eeq
Here the measure $Q$ is  absolutely continuous w.r.t. $P$, and
the likelihood process
$$
L_t=\frac{dQ}{dP} \quad\mbox{on $\mathcal{F}_t$},
$$
has dynamics  given by
$$
dL_t= L_t\varphi dW_t+ L_{t-}\int_{Z}\eta(z)\left[\mu(dz,dt)-\lambda(dz)dt\right],
$$
with $\varphi$ and $\eta$ given by
\beqarno
\varphi(t,s) &=& -\frac{\bracket{\alpha (t,s)-r+
\int_{Z}\beta(y) \lambda(dy)} \sigma (t,s)}{ [\sigma^{2}(t,s)+\int_{Z}\beta^2(y) \lambda(dy)]},\\
\eta (t,s,z)&=&- \frac{\bracket{\alpha (t,s)-r+ \int_{Z}\beta(y) \lambda(dy)} }
{[\sigma^{2}(t,s)+\int_{Z}\beta^2(y) \lambda(dy)]}\beta(z).
\eeqarno
From the Girsanov Theorem it follows that the $Q$ intensity $\lambda^Q$, of the point process $\mu(dt,dz)$ is given by
$$
\lambda^Q(t,s,dz)=\krull{1-\frac{\bracket{\alpha (t,s)-r+ \int_{Z}\beta(y) \lambda(dy)} }
{[\sigma^{2}(t,s)+\int_{Z}\beta^2(y) \lambda(dy)]}\beta(z)}\lambda (dz),
$$
and that
$$
dW_t=-\frac{\bracket{\alpha (t,S_t)-r+
\int_{Z}\beta(y) \lambda(dy)} \sigma (t,S_t)}{ [\sigma^{2}(t,S_t )+\int_{Z}\beta^2(y) \lambda(dy)]}dt +dW_t^Q
$$
where $W^Q$ is $Q$ a Wiener process.
A simple  calculation now shows that the $Q$ dynamics of the stock prices $S$ are given by
$$
dS_t=rS_tdt+S_t \sigma(t,S_t)dW_t^Q +S_{t-}\int_{\cal Z}\beta (z)\bracket{\mu (dt,dz)-\lambda^Q(t,S_t,dz)}
$$
so the measure $Q$ is in fact a risk neutral martingale measure,
and it is easy to check that $Q$ is in fact the so called ``minimal martingale measure''
used in the context of local risk minimization and developed in \cite{schweiz95} and related papers.
This fact was, in a Wiener process framework, observed already in \cite{Bas}.

Performing similar calculations, one can show that the solution of the PIDE \cref{Bpide1}-\cref{Bpide2}
can be represented as
\beqar
B(t,s)&=&\EtqX{t,s}{Q}{\int_{t}^{T} \frac{(\alpha-r+ \int_{Z}\beta(z) \lambda(dz))^2 }
{\gamma [\sigma^{2}+\int_{Z}\beta^2(z) \lambda(dz)]}d\tau}\nonumber
\\&&- \EtqX{t,s}{Q}{\int_{t}^{T} \half \sigma^2 \gamma [ e^{r(T-\tau)}\hat{u}+ b_{s}s]^2 d\tau} \nonumber \\&&
- \EtqX{t,s}{Q}{\int_{t}^{T}\frac{\gamma}{2}\int_{Z}[e^{r(T-\tau)}\hat{u}\beta(z)+\Delta_{\beta}b]^2\lambda(dz)d\tau},
\label{B_MPP2}
\eeqar
with $Q$ as above.

We can finally summarize our results.
\bprop
With notation as above, the following hold.
\bei
\item
The optimal amount of money invested in a stock is given by
$$
\hat{ \bf u} =\frac{(\alpha -r+ \int_{Z}\beta(z) \lambda(dz))}{\gamma [\sigma^{2}+\int_{Z}\beta^2(z) \lambda(dz)]}e^{-r(T-t)}-\frac{(\sigma b_{s}s+ \int_{Z}\beta(z) \Delta_{\beta}b(z)\lambda(dz))}
{ [\sigma^{2}+\int_{Z}\beta^2(z) \lambda(dz)]}e^{-r(T-t)}
$$
\item
The  mean-variance utility of the optimal portfolio is given by
$$
V(t,x) = e^{r(T-t)}x + b(t,s)
$$
where $b(t,s)$ is given by  stochastic representation \cref{b_MPP2}.
\item
The expected terminal value of the optimal portfolio  is
$$
\EnX{t,x}{X_T} = x e^{r(T-t)}+B(t,s)
$$
where $x$ is the present portfolio value and $B(t,s)$ is given by the stochastic representation  \cref{B_MPP2}.
\ei
\eprop
\section{Mean-variance control with wealth dependent risk aversion} \label{lmvr}
Going back to the MV portfolio optimization of Section \ref{mv},  there is a non trivial problem connected with the solution as presented in
Proposition \ref{mvprop}.
We recall that  the MV objective function at time $t$, given current wealth $X_t=x$, is given by
$$
\EnX{t,x}{X_T}-\frac{\gamma}{2} \VnX{t,x}{X_T},
$$
where $X_T$ is the wealth at the end of the time period, and where $\gamma$ is a given  constant
representing the risk aversion of the agent. Referring to Proposition \ref{mvprop} we see that, at time $t$ and when total wealth is $x$, the optimal dollar amount $u(t,x)$ invested in the risky asset is of the form
$$
u(t,x)=h(t)
$$
where $h$ is a deterministic function of time. In particular this implies that the dollar   amount  invested  in the risky asset
does {\bf not} depend on current wealth $x$.

This result is economically unreasonable, and stems from the fact that the risk aversion parameter $\gamma$ above is assumed to be a constant.
A more natural idea is  that we should explicitly {\bf allow  $\mathbf \gamma$
 to depend on current wealth}.

We thus go on  to study the MV problem
with a state dependent risk aversion.
More explicitly we consider an objective function of the form
$$
\EnX{t,x}{X_T}-\frac{\gamma (x)}{2} \VnX{t,x}{X_T},
$$
where $\gamma$ is a deterministic function of present wealth $x$.  The financial market is as in Section \ref{mvs} so the wealth dynamics are again given by
$$
dX_t= [rX_t+(\alpha-r)u_t]dt+\sigma u_tdW_t.
$$
This type of problem cannot easily be treated within the framework
of \cite{Bas} or \cite{Czi}, but it is a simple special case of the theory developed in the present paper.
We now go on to treat this improved model formulation for mean variance portfolio selection, and we follow our previous paper \cite{bmz}.
See also \cite{bjo-murFSD} for a discrete time version.
\subsection{The extended HJB equation for a general $\gamma (x)$} \label{gam}
Applying the general theory from Section \ref{cont} to  the objective  functional
\beqno
J(t,x, \op{u})=\EnX{t,x}{X_T^{\op{u}}}-\frac{\gamma (x)}{2}\VnX{t,x}{X_T^{\op{u}}}
\eeqno
gives us, after a large number of elementary calculations,  the following result (see \cite{bmz} for details).
We use the notation $\beta = \alpha -r$.

\bprop \label{nprop1}
The extended HJB system takes the following form:
\beqarno
V_t + \sup_{u \in {R}}
\left \{ \pa{rx+ \beta u}\bracket{ V_x- f_y-\frac{\gamma'(x)}{2}g^2}    \right.  \nonumber\\
\left.  + \half   \sigma^2 u^2  \bracket{V_{xx} - f_{yy} - 2 f_{xy}
-\frac{\gamma''(x)}{2}g^2-2\gamma '(x)gg_x  -\gamma (x)g_x^2  }  \right \} & = & 0, \label{hjb4}\\
f_t(t,x,y)+ \pa{rx+ \beta \hat{u}} f_x(t,x,y) +\half {\sigma}^2 \hat{u}^2f_{xx}(t,x,y)&=&0, \label{f4} \\
g_t(t,x)+\pa{rx+ \beta \hat{u}} g_x(t,x) +\half {\sigma}^2 \hat{u}^2 g_{xx}(t,x)&=&0. \label{g4}
\eeqarno
\eprop

\subsection{A special choice of $\gamma (x)$} \label{gams}
One way of finding a natural candidate of $\gamma$ is to perform a dimension analysis. In the reward function
\beqno
J(t,x, \op{u})=\EnX{t,x}{X_T^{\op{u}}}-\frac{\gamma (x)}{2}\VnX{t,x}{X_T^{\op{u}}}
\eeqno
we see that the first term
$
\EnX{t,x}{X_T^{\op{u}}}
$
has the dimension $(dollar)$. The term $\VnX{t,x}{X_T^{\op{u}}}$ has the dimension $(dollar)^2$, so in order to have
a reward function measured in dollars we have to choose $\gamma$ in such a way that  $\gamma (x)$ has the dimension $(dollar)^{-1}$.
The most obvious way to accomplish this is of course to specify $\gamma$ as
\beqno
\gamma (x)=\frac{\gamma}{x},
\eeqno
where, with a slight abuse of notation, the occurrence of $\gamma$ in the right hand side denotes a constant. This is also very natural
from an economic point of view since the risk aversion will now be decreasing in wealth.

We now conjecture that $\hat{u}$ is linear in $x$ so we make the {\em Ansatz}
$$
\hat{u}(t,x)= c(t)x
$$
for some deterministic function $c$.
If this is the case, then  $X$  will be  GBM so we will have
\beqar
\EnX{t,x}{X_T^{\hat{\op{u}}}}&=&a(t)x, \label{s4}\\
\EnX{t,x}{\pa{X_T^{\hat{\op{u}}}}^2}&=&b(t)x^2, \label{s5}
\eeqar
for some deterministic functions $a$ and $b$.
We recall the probabilistic interpretations
\beqarno
f(t,x,y)&=&\EnX{t,x}{X_T^{\hat{\op{u}}}} -\frac{\gamma (y)}{2}\EnX{t,x}{\pa{X_T^{\hat{\op{u}}}}^2}, \label{s6}\\
g(t,x)&=&\EnX{t,x}{X_T^{\hat{\op{u}}}}.\label{s7}
\eeqarno
This leads to the {\em Ansatz}
\beqarno
f(t,x,y)&=&a(t)x -\frac{\gamma }{2y}b(t)x^2, \label{s8}\\
g(t,x)&=&a(t)x.\label{s9}
\eeqarno

Plugging these expressions into the extended HJB system gives us the following results.

\bprop \label{natprop}
For the case when $\gamma (x)=\frac{\gamma}{x}$, the equilibrium control is given by
$$
\hat{u}(t,x)=\frac{\beta}{\gamma \sigma^2 }\frac{a(t) + \gamma \bracket{a^2(t) -b(t)}}{b(t)}x
$$
where $a$ and $b$ solves the ODE system
\beqarno
\dot{a}+\pa{r+  \frac{\beta^2}{\gamma \sigma^2 b}\pa{a + \gamma \bracket{a^2 -b}}}a&=&0,\\
a(T,x)&=&1,\\
\dot{b} +\krull{2\pa{r+  \frac{\beta^2}{\gamma \sigma^2 b}\pa{a + \gamma \bracket{a^2 -b}}}
+ \frac{\beta^2}{\gamma^2 \sigma^2 b^2}\pa{a + \gamma \bracket{a^2 -b}}^2}b&=&0,\\
b(T,x)&=&1.
\eeqarno
The equilibrium value function $V$ is given by
$$
V(t,x)=\krull{a(t)+ \frac{\gamma}{2}\bracket{a^2(t)-b(t)}}x.
$$
\eprop
\begin{remark}
We see that the optimal dollar value $\hat{u}$ invested in the risky asset is now linear in wealth.
From an economic point of view this is a much nicer result than the one obtained from the `´naive'' mean variance model where
$\hat{u}$ was constant as a function of $x$.
\end{remark}

\section{The  inconsistent linear quadratic regulator} \label{lqr}
To illustrate how the theory works in another simple case, we consider a small variation of the classical linear quadratic regulator.
The model is specified as follows.
\bei
\item
The value functional for Player $t$ is given by
$$
\EnX{t,x}{\half \int_t^Tu_s^2ds} + \frac{\gamma}{2}\EnX{t,x}{\pa{X_T-x}^2}
$$
where $\gamma$ is a positive constant.
\item
The state process $X$ is  scalar with dynamics
$$
dX_t=\bracket{aX_t + bu_t}dt+\sigma dW_t,
$$
where $a$, $b$ and $\sigma$ are given constants.
\item
The control $u$ is scalar with no constraints.
\ei
This is a time inconsistent version of the classical linear quadratic regulator.
Loosely speaking we want control the system so that the final sate $X_T$ is close to $x$ while at the same time
keeping the control energy (formalized by the integral term) is small. The time inconsistency stems from the fact that
the target point $x=X_t$ is changing as time goes by. In discrete time this problem is studied in \cite{bjo-murFSD}.

For this problem we have
\beqarno
F(x,y)&=&{\frac{\gamma}{2}(y-x)^2},\\
{H}(u)&=&\half u^2
\eeqarno
and as usual we introduce  the functions $f^y(t,x)$ and $f(t,x,x)$ by
\beqarno
f^y(t,x)&=&{\EnX{t,x}{\int_t^T\frac{1}{2}\hat{\op{u}}^2_s\left(X_s^{\hat{\op{u}}} \right)ds+\pa{X_T^{\hat{\op{u}}}-y}^2}},\\
f(t,x,y)&=&f^y(t,x).
\eeqarno
The extended HJB system takes the form
\beqarno
{\inf_u \krull{ \half u^2 + \op{A}^uf^x(t,x)}}&=&0, \quad 0 \leq t \leq T \\
\op{A}^{\hat{\op{u}}}f^y(t,x)+{\frac{1}{2}\hat{\op{u}}^2_t(x)}&=&0, \quad 0 \leq t \leq T \\
f^y(T,x)&=&\frac{\gamma}{2}(x-y)^2.
\eeqarno
From the $X$ dynamics we see that
$$
\op{A}^u=\dfdx{}{t} + (ax+bu)\dfdx{}{x}+ \half \sigma^2 \ddfdxx{}{x}.
$$
We thus obtain the following form of the HJB equation, where for shortness of notation we denote partial derivatives
by lower case index so, for example, $f_x=\dfdx{f}{x}$. 
\beqarno
{\inf_u \krull{ \half u^2 +f_t(t,x,x) +(ax+bu)f_x(t,x,x) +\frac{1}{2}\sigma^2 f_{xx}(t,x,x)}}&=&0,\\
{f(T,x,x)}&=&0.
\eeqarno
The coupled system for $f^y$ is given by
\beqarno
f_t^y(t,x)+ \bracket{ax + b \hat{\op{u}}(t,x)}f_x^y(t,x) + \half \sigma^2 f_{xx}^y(t,x)+{\frac{1}{2}\hat{\op{u}}^2(t,x)}&=&0,\\
f^y(T,x)&=&\frac{\gamma}{2}(x-y)^2.
\eeqarno
The first order condition in the HJB equation gives us
$$
\hat{\op{u}}(t,x)={-b f_x(t,x,x)}
$$
and, inspired of the standard regulator problem,  we now make the Ansatz (attempted solution)
\beqar
f(t,x,y)&=&A(t)x^2 + B(t)y^2 + C(t)xy + D(t)x + F(t)y+H(t), \label{lq1}
\eeqar
where all coefficients are deterministic functions of time.
We now insert the Ansatz into the HJB system, and perform a number of extremely boring calculations.
As a result of these calculations, it turns out that the variables separate in the expected way and we have the following result.
\bprop
For the time inconsistent regulator, we have the structure \cref{lq1}, where the coefficient functions solve the following system of ODE:s.
\beqarno
A_t+2aA-2b^2A(2A+C)+\frac{1}{2}b^2(2A+C)^2&=&0\\
B_t&=&0\\
C_t+aC-b^2C(2A+C)&=&0\\
D_t+aD-2b^2AD &=&0 \\
F_t-b^2CD&=&0\\
H_t-\frac{1}{2}b^2D^2 
+\sigma^2 A&=&0
\eeqarno
With boundary conditions
\beqarno
A(T)=\frac{\gamma}{2}, \quad B(T)&=&\frac{\gamma}{2}, \quad C(T)=-\gamma,\\
D(T)=0,\quad  F(T)&=&0, \quad H(T)=0.
\eeqarno
\eprop

\section{A Cox-Ingersoll-Ross production economy with time inconsistent preferences} \label{eqp}
In this section we apply the previously developed theory  to a rather detailed study of a general equilibrium model for a production economy with time inconsistent preferences. The model under consideration is a time inconsistent analogue of the classic Cox-Ingersoll-Ross model in \cite{CIR85a}.  Our main objective is to investigate the structure of the equilibrium short rate, the equilibrium Girsanov kernel, and the equilibrium stochastic discount factor. 

There are a few earlier papers on equilibrium with time inconsistent preferences. In \cite{Bar99} and \cite{lutt03} the authors study continuous time equilibrium models of a particular type of time inconsistency, namely non-exponential discounting. While \cite{Bar99} considers a deterministic neoclassical model of economic growth, \cite{lutt03} analyze general equilibrium in a stochastic endowment economy. 

Our present study is much inspired by the earlier paper  \cite{Khapko15} which, in very great detail,  studies equilibrium in a very general setting of an endowment economy with dynamically inconsistent preferences that is not limited to the particular case of non-exponential discounting. 

Unlike the  papers mentioned above, which all studies {\em endowment} models, 
we study a stochastic \emph{production} economy of Cox-Ingersoll-Ross type.  

\nocite{Harris12}
\subsection{The  Model} \label{eqm}
We start with some formal assumptions concerning the production technology.

\begin{assumption}
We assume that  there exists a constant returns to scale physical production technology process $S$ with dynamics
\beq \label{seq}
dS_t=\alpha S_tdt + S_t\sigma dW_t.
\eeq
The economic agents can invest unlimited positive amounts in this technology, but since it is a matter of physical investment, short positions are not allowed. 
\end{assumption}

More concretely this means that at any time you are allowed to invest dollars in the production process. 
If you, at time $t$, invest $q$ dollars, and wait until time $u$ then you will receive the amount of
$q \cdot {S_{u}}/{S_{t}}$
in dollars. In particular we see that the return on the investment is linear in $q$, hence the term ``constant returns to scale''.
Since this is a matter of physical investment, short-selling is not allowed.

A moment of reflection  shows that, from a purely formal point of view, investment in the technology $S$ 
is in fact  equivalent to the possibility of investing in a risky asset with price process $S$, but again 
with the constraint that short-selling is not allowed.

We also need a risk free asset, and this is provided by the next assumption.
\begin{assumption}
We  assume that there exists a risk free asset in zero net supply with dynamics
$$
dB_t=r_tB_tdt,
$$
where $r$ is the short rate process, which will be determined endogenously. The risk free rate $r$ is assumed to be of the form
$$
r_t=r(t,X_t)
$$
where $X$ denotes the portfolio value of the representative investor (to be defined below).
\end{assumption}

Interpreting the production technology $S$ as above,  the wealth dynamics will be given by
$$
dX_t=X_tu_t(\alpha -r_t)dt + (r_tX_t-c_t)dt + X_tu_t\sigma dW_t.
$$
where $u$ is the portfolio weight on production, so $1-u$ is the weight on the risk free asset.
Finally we need an economic agent.

\begin{assumption}
We  assume that there exists a representative agent 
who, at every point $(t,x)$,  wishes to maximize the reward functional
\beq \label{eqm10}
\EtqX{t,x}{}{\int_t^T U(t,x,s, c_s)ds}.
\eeq
\end{assumption}

\subsection{Equilibrium definitions} \label{edp}
We now go on to study equilibrium in our model.  We will in fact have two equilibrium concepts
\bei
\item
Intrapersonal equilibrium.
\item 
Market equilibrium.
\ei
The intrapersonal equilibrium is related to the lack of time consistency in preferences, 
whereas the market equilibrium is related to market clearing. We now discuss these concepts in more detail.
\subsubsection{Intrapersonal equilibrium} \label{ipe}
Consider, for a given short rate function $r(t,x)$ the control problem with reward functional
$$
\EtqX{t,x}{}{\int_t^T U(t,x,s,c_s)ds}.
$$
and wealth dynamics
$$
dX_t=X_tu_t(\alpha -r_t)dt + (r_tX_t-c_t)dt + X_tu_t\sigma dW_t.
$$
where $r_t$ is shorthand for $r(t,X_t)$.
If the agent wants to maximize the reward functional for every initial point $(t,x)$ then, because of the appearance of $(t,x)$ in the utility function $U$, this is a time inconsistent control problem.
In order to handle this situation we use the game theoretic setup and results  developed in Sections \ref{int}-\ref{cex} above. This subgame perfect Nash equilibrium concept   is henceforth referred to as the {\bf intrapersonal equilibrium}. 
\subsubsection{Market equilibrium} \label{med}
By a {\bf market equilibrium} we  mean a situation where the agent follows an intrapersonal equilibrium strategy, and where  the market clears for the risk free asset. The formal definition is as follows.
\bdf
A market  equilibrium of the model is a  triple of real valued functions $\krull{\hat{c}(t,x),\hat{u}(t,x),r(t,x)}$  such that the following hold.
\ben
\item
Given the risk free rate process $r(t,x)$, the intrapersonal equilibrium consumption and investment are given by $\hat{c}$ and $\hat{u}$ respectively.
\item 
The market clears for the risk free asset, i.e.
$$
\hat{u}(t,x)\equiv 1.
$$
\en
\edf
\subsection{Main goals of the study} \label{map}
As will be seen below, there will be a unique equilibrium martingale measure $Q$ with corresponding likelihood process $L =dQ/dP$, where $L$ has dynamics
$$
dL_t=L_t\vp_tdW_t.
$$
The process $\vp$ will be referred to as the equilibrium {\em Girsanov kernel}. There will also be a  equilibrium short rate process $r$, which will be related to $\vp$ by the standard no arbitrage relation 
\beq \label{map1}
r(t,x)= \alpha + \vp (t,x)\sigma ,
\eeq
which says that $S/B$ is a $Q$-martingale.
There will also be a unique equilibrium stochastic discount factor $M$ defined by
$$
M_t=e^{-\int_0^tr_sds}L_t.
$$
For ease of notation we will, however, only identify the stochastic discount factor $M$,  up to a multiplicative constant, 
so for any arbitrage free (non dividend) price process $p_t$ we will have the pricing equation
$$
p_s=\frac{1}{M_s}\EpXF{P}{M_tp_t}{\Ft{s}}.
$$  
Our goal is to obtain expressions for $\vp$, $r$ and $M$. 
\subsection{The extended HJB equation} \label{eh}
In order to determine the intrapersonal equilibrium we use  the results from Section \ref{GO}.  
The extended HJB equation \cref{cexv11} now reads as
\beq \label{eh1}
 \sup_{u {\geq0}, c \geq 0}   \krull{  U(t,x,t,c)  + \op{A}^{u,c}f^{tx}(t,x)} = 0
\eeq
and $f^{sy}$ is determined by 
\beq \label{eh31}
\op{A}^{\hat{\op{u}},\hat{\op{c}}}f^{sy}(t,x)+U\pa{s,y,t,\hat{c}(t,x)}=0
\eeq
with the probabilistic representation 
\beq \label{eh3}
f^{sy}(t,x)=\EnX{t,x}{\int_t^TU\pa{s,y,\tau,\hat{\op{c}}(\tau, X_{\tau}^{\hat{\op{u}}})}d{\tau}}, \quad 0 \leq t \leq T.
\eeq

The term $\op{A}^{u,c}f^{tx}(t,x)$ is given by
\beq \label{eh7}
\op{A}^{u,c}f^{tx}(t,x)=f_t+ xu\pa{\alpha -r}f_x+ (rx-c)f_x+\half x^2u^2 \sigma^2 f_{xx}
\eeq
where $f$ and the derivatives are evaluated at $(t,x,t,x)$ and where we have used the notation
\beqarno
f(t,x,s,y)&=&f^{sy}(t,x),\\
f_t(t,x,s,y)&=&\dfdx{f}{t}(t,x,s,y),\\
f_x(t,x,s,y)&=&\dfdx{f}{x}(t,x,s,y),\\
f_{xx}(t,x,s,y)&=&\ddfdxx{f}{x}(t,x,s,y).
\eeqarno
The first order conditions for an interior optimum are
\beqar
U_c(t,x,t,\hat{c})&=&f_x(t,x,t,x) \label{eh10}\\
\hat{u}&=&-\pa{\frac{\alpha -r}{\sigma^2}}\frac{f_x(t,x,t,x)}{xf_{xx}(t,x,t,x)}.\label{eh11}
\eeqar

\subsection{Determining market equilibrium} \label{me}
In order to determine the market equilibrium we use the equilibrium condition $\hat{u}=1$. Plugging this into \cref{eh11} we immediately obtain 
our first result.
\bprop \label{me1}
With assumptions as above the following hold.
\bei
\item
The equilibrium short rate is given by
\beq \label{me2}
r(t,x)=\alpha + \sigma^2 \frac{xf_{xx}(t,x,t,x)}{f_{x}(t,x,t,x)}.
\eeq
\item
The equilibrium Girsanov kernel $\vp$ is given by
\beq \label{me3}
\vp(t,x)=\sigma \frac{xf_{xx}(t,x,t,x)}{f_{x}(t,x,t,x)}.
\eeq
\item
The extended equilibrium HJB system has the form
\beqar \label{me4}
U(t,x,t,\hat{c}) +f_t+ \pa{\alpha x -\hat{c}}f_x+\half x^2 \sigma^2 f_{xx}&=&0,\\
&& \nonumber\\
\op{A}^{\hat{\op{c}}}f^{sy}(t,x)+U\pa{s,y,t,\hat{c}(t,x)}&=&0
\eeqar
\item
The equilibrium consumption $\hat{c}$ is determined by the first order condition 
\beq
U_c(t,x,t,\hat{c})=f_x(t,x,t,x)
\eeq
\item
The term $\op{A}^{\hat{\op{c}}}f^{tx}(t,x)$ is given by
\beq 
\op{A}^{\hat{\op{c}}}f^{tx}(t,x)=f_t+ x\pa{\alpha -r}f_x+ (rx-\hat{c})f_x+\half x^2 \sigma^2 f_{xx}
\eeq
\item
The equilibrium $X$ dynamics are given by
\beq \label{me5}
dX_t=\pa{\alpha X_t -\hat{c}_t}dt  + X_t\sigma dW_t.
\eeq
\ei
\eprop

\proof 
The formula \cref{me3} follows from \cref{me2} and \cref{map1}. The other results are obvious. \endproof
\subsection{Recap of standard results} \label{rsdf}
We can  compare the results above with the standard case where the utility functional for the agent  is of the time consistent form 
$$
\EtqX{t,x}{}{\int_t^T U(s, c_s)ds}.
$$
In this case we  have a standard HJB equation of the form
\beq
 \sup_{u\in R, c \geq 0}   \krull{  U(t,c)  + \op{A}^{u,c}V(t,x)} = 0.
\eeq
and the equilibrium quantities are given by the well known expressions
\beqar
r(t,x)&=&\alpha + \sigma^2 \frac{xV_{xx}(t,x)}{V_{x}(t,x)},\label{me10}\\
\vp(t,x)&=&\sigma \frac{xV_{xx}(t,x)}{V_{x}(t,x)}.\label{me11}
\eeqar
We note the strong structural similarities between the old and the new formulas, but we also note important differences. 
Let us take the formulas for the equilibrium short rate $r$ as an example. We recall the standard and time inconsistent formulas
\beqar
r(t,x)&=&\alpha + \sigma^2 \frac{xV_{xx}(t,x)}{V_{x}(t,x)},\label{rsdf20}\\
r(t,x)&=&\alpha + \sigma^2 \frac{xf_{xx}(t,x,t,x)}{f_{x}(t,x,t,x)}.\label{rsdf21}
\eeqar
For the time inconsistent case we do have the relation $V^e(t,x)=f(t,x,t,x)$ (where temporarily, and for the sake of clearness,
 $V^e$ denotes the equilibrium value function)
so it is tempting to think that we should be able to write \cref{rsdf21} as
$$
r(t,x)=\alpha + \sigma^2 \frac{xV_{xx}^e(t,x)}{V_{x}^e(t,x)}
$$
which would be structurally identical to \cref{rsdf20}. Not however, that while we do have $f(t,x,t,x)=V^e(t,x)$,
the partial derivative $f_x(t,x,t,x)$ is {\bf not} equal to $V_x^e(t,x)$. The reason is that  in $f_x(t,x,t,x)$, the partial differential operator only acts on the first occurrence of $x$ in $f(t,x,t,x)$, whereas $V_x^e(t,x)$ will 
be the total $x$-derivative of $f(t,x,t,x)$.  
\subsection{The stochastic discount factor} \label{sdf}
We now go on to investigate our main object of interest, namely the equilibrium stochastic discount factor $M$.
We recall from general arbitrage theory that
\beq \label{sdf1}
M_t=e^{-\int_0^tr_udu}L_t
\eeq
where $L$ is the likelihood process
$$
L_t=\frac{dQ}{dP}, \quad \mbox{on $\Ft{t}$}
$$
with dynamics
$$
dL_t=L_t\vp_tdW_t.
$$
From this we immediately obtain the $M$ dynamics as
\beq \label{sdf5}
dM_t=-r_tM_tdt + M_t\vp_tdW_t,
\eeq
so we can identify the short rate $r$ and the Girsanov kernel $\vp$ from the dynamics of $M$. 

From Proposition \ref{me1} we know $r$ and $\vp$, so in principle we have in fact already determined $M$, but 
we now want to investigate the relation between $M$, the direct utility function $U$, and   the indirect utility function $f$ in the extended HJB equation. 

We  recall from standard theory that for the usual time consistent case the (non normalized) stochastic discount factor $M$ is given by
\beqno
M_t=V_x(t,X_t),
\eeqno
or equivalently by
\beqno
M_t=U_c(t,c_t)
\eeqno
along the equilibrium path.
In our present setting we have
$$
V(t,x)=f(t,x,t,x),
$$
so a conjecture would perhaps be that the  stochastic discount factor  for the time inconsistent case is given by
at least one  of the formulas
\beqarno
M_t&=&V_x(t,X_t),\\
M_t&=&f_x(t,X_t,t,X_t),\\
M_t&=&U_c(t,X_t,t,c_t)
\eeqarno
along the equilibrium path. In order to check if any of these formulas are correct we only have to compute the 
corresponding differential $dM_t$ and check whether it satisfies 
\cref{sdf5}. It is then easily seen that none of the formulas for $M$ are correct.   The situation is in thus  more complicated and we now go on to derive the correct representation of the stochastic discount factor.
\subsubsection{A representation formula for $M$} \label{m}
We now go back to the time inconsistent case with  utility of the form
$$
\EtqX{t,x}{}{\int_t^T U(t,x,s,c_s)ds}.
$$
We will, below, present a representation for the stochastic discount factor $M$ in the production economy, but first we need to introduce some new notation.

\bdf \label{m2}
Let $X$ be a (possibly vector valued) semimartingale and let $Y$ be an optional process. For a $C^{2}$ function $f(x,y)$ we introduce the 
``partial stochastic differential'' $\partial_x$ by the formula
\beq \label{m4}
\partial_x f(X_t,Y_t)=df(X_t,y),\quad \mbox{evaluated at $y=Y_t$.}
\eeq
\edf
The intuitive interpretation of this is  that
\beq \label{m41}
\partial_x f(X_t,Y_t)=f(X_{t+dt},Y_t)-f(X_t,Y_t),
\eeq
and we have the following proposition, which generalizes the standard result for the time consistent theory.

\bth \label{m1}
The stochastic discount factor $M$ is determined by
\beq \label{m7}
d\pa{\ln M_t}=\partial_{t,x}\ln f_x(t,X_t,t,X_t),
\eeq
where the partial differential $\partial_{t,x}$ only operates on the variables $(t,x)$ in $f_x(t,x,s,y)$.
Alternatively we can write
\beq \label{m71}
M_t=U_c\pa{t,X_t,t,\hat{c}_t} \cdot e^{Z_t},
\eeq
where $Z$ is determined by
\beq \label{m72}
dZ_t=\partial_{tx}\ln f_x \pa{t,X_t,t, X_t}- d\ln f_x \pa{t,X_t,t,X_t}.
\eeq
\eth

\begin{remark}
We remark here on the structural similarity of the stochastic discount factor to the result obtained in \cite{Khapko15}. \end{remark}

\begin{remark}
For a more concrete  interpretation of this result, see Section \ref{mi} below.
Note again that the operator $\partial_{tx}$ in \cref{m72} only acts on the first occurrence  of $t$ and $X_t$ in  
$f_x \pa{t,X_t,t, X_t}$ whereas 
the operator $d$ acts on  the entire process $t \longmapsto f_x\pa{t,X_t,t, X_t}$.
\end{remark}

\proof
Formulas \cref{m71}-\cref{m72} follows from \cref{m7} and the first order condition 
$U_c\pa{t,X_t,t,\hat{c}_t} =f_x \pa{t,X_t,t, X_t}$. It thus remains to prove \cref{m7}.

From \cref{sdf5}  it follows that we need to show that
\beq \label{m10}
\partial_{t,x}\ln f_x(t,X_t,t,X_t)=-\krull{r_t + \half \vp_t^2}dt + \vp_t dW_t
\eeq
where $r$ and $\vp$ are given by \cref{me2}-\cref{me3}.
Applying Ito and the definition of $\partial_{t,x}$  we obtain
$$
\partial_{t,x}\ln f_x(t,X_t,t,X_t)=A(t,X_t)dt + B(t,X_t)dW_t,
$$
where
\beqar
A(t,x)&=&\frac{1}{f_x}\krull{f_{xt} + \pa{\alpha x - \hat{c}}f_{xx}+ \half \sigma^2 x^2 f_{xxx} - \half \sigma^2 x^2 \frac{f_{xx}^2}{f_x}}, \label{m11}\\
B(t,x)&=&\sigma x \frac{f_{xx}}{f_x}. \label{m12}
\eeqar
From \cref{me3} we see that indeed $B(t,x)=\vp (t,x)$ so, using \cref{me2},  it remains to show that
\beq \label{m13}
A(t,x)=-\krull{\alpha + \sigma^2 x\frac{f_{xx}}{f_x} + \half \sigma^2 x^2  \frac{f_{xx}}{f_x} }.
\eeq
To show this we differentiate the equilibrium HJB equation \cref{me4}, use the first order condition $U_c=f_x$,  and obtain
\beqar
U_x + f_{ty}+f_{tx}+ \pa{\alpha x - \hat{c}}f_{xx} + \pa{\alpha x - \hat{c}}f_{xy} +\alpha f_x && \nonumber\\
+ \sigma^2 x f_{xx} + \half \sigma^2 x^2 f_{xxy} + \half \sigma^2 x^2 f_{xxx}&=&0, \label{m15}
\eeqar
where $f_{tx}=f_{tx}(t,x,t,x)$ and similarly for other derivatives,  $\hat{c}=\hat{c}(t,x)$ and $U_x=U_x\pa{t,x,t,\hat{c}(t,x)}$.
From the extended HJB system we also recall the PDE for $f^{sy}$ as
$$
f_t^{sy}(t,x) + \pa{\alpha x - \hat{c}}f_x^{sy}(t,x) + \half \sigma^2 x^2 f_{xx}^{sy}(t,x)+ U(s,y,t,\hat{c})=0
$$
Differentiating this equation w.r.t. the variable $y$ and evaluating at $(t,x,t,x)$ and $\hat{c}(t,x)$ we obtain
\beq \label{m16}
f_{ty} + \pa{ \alpha x-\hat{c}}f_{xy} + \half \sigma^2 x^2 f_{xxy} + U_x=0.
\eeq
We can now plug this into \cref{m15} to obtain
$$
f_{tx}+ \pa{ \alpha x -\hat{c}}f_{xx} + \alpha f_x + \sigma^2 x f_{xx} + \half \sigma^2 x^2 f_{xxx}=0.
$$ 
Plugging this into \cref{m11} we can write $A$ as
$$
A(t,x)=- \krull{\alpha  + \sigma^2 x \frac{f_{xx}}{f_x} + \half \sigma^2 x^2 \frac{f_{xx}^2}{f_x}}
$$
which is exactly \cref{m13}. \endproof
\subsubsection{Interpreting the representation formula} \label{mi}
The representation formula \cref{m7} does not {\em prima facie}  seem to have a natural interpretation. We can of course write $M$  as
$$
M_t=e^{\int_0^t \partial_{t,x} \ln f_x(s,X_s,s,X_s)}
$$
but this does not seem to give much insight. In order to get a deeper understanding we recall that for any (non dividend) asset price 
process $p$ we have the valuation formula
$$
p_s=\EpXF{P}{\frac{M_t}{M_s}p_t}{\Ft{t}}.
$$
it is therefore natural to make the following definition.

\bdf
For any $s < t$ we define the $(s,t)$-stochastic  discount factor $M_{st}$ by
\beq
M_{st}=\frac{M_t}{M_s}.
\eeq
\edf 
We thus have a natural multiplicative  structure in the sense that for $s < u < t$ we have
$$
M_{st}=M_{su}\cdot M_{ut}.
$$
The economic interpretation of $M_{st}$ is thus that (via conditional expectation) it discounts the value of a 
stochastic claim at time $t$ back to time $s$. 
It is now natural to look at the infinitesimal version of $M_{st}$. This object
would intuitively be defined by the formula
\beq \label{mi5}
M_{t,t+dt}=\frac{M_{t+dt}}{M_t}
\eeq
and it would tell us how we discount on the infinitesimal scale from time $t+dt$ back to time $t$. 

In order to make  \cref{mi5} more precise we note that we can write it as
$$
M_{t,t+dt}=e^{\ln M_{t+dt}-\ln M_t}=e^{d\pa{\ln M_t}}
$$
and this motivates the following formal definition.
\bdf
The  log stochastic discount factor $m_t$ is defined by
\beq \label{mi10}
m_t={\ln M_t} 
\eeq
\edf
We thus have 
$$
M_{st}=e^{m_t-m_s}
$$
and the informal interpretation
$$
M_{t,t+dt}=e^{dm_t}.
$$
Theorem \ref{m1} shows that
$$
dm_t=\partial_{t,x}\ln f_x(t,X_t,t,X_t),
$$
so we have the interpretation
$$
M_{t,t+dt}=e^{\partial_{t,x}\ln f_x(t,X_t,t,X_t)}
$$
Using the interpretation \cref{m41} and doing some simple calculations we finally obtain the following (informal) result.
\bprop \label{mi12}
With notation as above we  have the following informal representation, corresponding to equation \cref{m7}.
\beq \label{mi15}
\frac{M_{t+dt}}{M_t}=\frac{f_x(t+dt,X_{t+dt},t,X_t)}{f_x(t,X_t,t,X_t)}.
\eeq
This can also be written as
\beq \label{mi17}
\frac{M_{t+dt}}{M_t}=\frac{U_c(t+dt,t+dt,X_{t+dt},\hat{c}_{t+dt})}{U_c(t,t,X_{t},\hat{c}_{t})}\cdot 
\frac{f_x(t+dt,X_{t+dt},t,X_t)}{f_x(t+dt,X_{t+dt},t+dt,X_{t+dt})}
\eeq
corresponding to equations \cref{m71}-\cref{m72}.
\eprop
Formula \cref{mi15} has a natural economic interpretation which can be seen from a dimension argument. 
The valuation formula for a price process $p$ will, on  the infinitesimal scale, read as
$$
p_t=\EpXF{P}{\frac{f_x(t+dt,X_{t+dt},t,X_t)}{f_x(t,X_t,t,X_t)}p_{t+dt}}{\Ft{t}}
$$
If we denote the dimension of money by dollars we have the following relations, where $dim$ denotes dimension
\beqarno
dim\bracket{ f_x(t+dt,X_{t+dt},t,X_t)}&=&\mbox{marginal utility, at $t$, of dollars at $t+dt$.}\\
dim \bracket{f_x(t,X_t,t,X_{t})}&=&\mbox{marginal utility, at $t$, of dollars at $t$.}
\eeqarno 
Since $p_{t+dt}$ has dimension dollars at $t+dt$, we see that by multiplying with the  factor $f_x(t,X_t,t+dt,X_{t+dt})$ 
transforms this dollar amount into marginal utility at time $t$. Dividing by $f_x(t,X_t,t,X_{t})$ gives us dollars at $t$, which is the dimension of $p_t$.
\subsection{Equilibrium with non-exponential discounting} \label{eqned}
A  case of particular interest occurs when the utility function is of the form
$$
U(t,x,s,c_s)=\beta (s-t)U(c_s)
$$
so the utility functional has the form
$$
\EtqX{t,x}{}{\int_t^T \beta (s-t)U(c_s)ds}.
$$

\subsubsection{Generalities} \label{neg}
 In the case of non exponential discounting it is natural to consider 
the case with infinite horizon. We will thus assume that $T= \infty$ so we have the functional
\beq \label{ned1}
\EtqX{t,x}{}{\int_t^\infty \beta (\tau-t)U(c_{\tau})d\tau}.
\eeq
The function $f(t,x,s,y)$ will now be of the form $f(t,x,s)$ and, 
because of the time invariance,  it is  natural to look for time invariant equilibria where
\beqarno
\hat{\op{u}}(t,x)&=&\hat{\op{u}}(x),\\
V(t,x)&=&V(x),\\
f(t,x,s)&=&g(t-s,x),\\
V(x)&=&g(0,x).
\eeqarno
Observing that $f_x(t,x,t)=g_x(0,x)=V_x(x)$ and similarly for second order derivatives, we may now restate proposition \ref{me1}.

\bprop \label{ned5}
With assumptions as above the following hold.
\bei
\item
The equilibrium short rate is given by
\beq \label{ned6}
r(x)=\alpha + \sigma^2 \frac{xV_{xx}(x)}{V_{x}(x)}.
\eeq
\item
The equilibrium Girsanov kernel $\vp$ is given by
\beq \label{ned7}
\vp(x)=\sigma \frac{xV_{xx}(x)}{V_{x}(x)}.
\eeq
\item
The extended equilibrium HJB system has the form
\beqar \label{ned8}
U(\hat{c}) +g_t(0,x)+ \pa{\alpha x -\hat{c}}g_x(0,x)+\half x^2 \sigma^2 g_{xx}(0,x)&=&0,\\
&& \nonumber\\
\op{A}^{\hat{\op{c}}}g(t,x)+\beta (t)U\pa{\hat{c}(x)}&=&0, \label{ned81}
\eeqar
\item
The function $g$ has the representation
\beq \label{ned9}
g(t,x)=\EnX{0,x}{\int_0^{\infty}\beta (t+ s)U\pa{\hat{c}_s}ds}
\eeq
\item
The equilibrium consumption $\hat{c}$ is determined by the first order condition 
\beq \label{ned10}
U_c(\hat{c})=g_x(0,x)
\eeq
\item
The term $\op{A}^{\hat{\op{c}}}g(t,x)$ is given by
\beq 
\op{A}^{\hat{\op{c}}}g(t,x)=g_t(t,x)+ x\bracket{\alpha -\hat{c}(x)}g_x(t,x)+ \half x^2 \sigma^2 g_{xx}(t,x)
\eeq
\item
The equilibrium $X$ dynamics are given by
\beq \label{ned11}
dX_t=\pa{\alpha X_t -\hat{c}_t}dt  + X_t\sigma dW_t.
\eeq
\ei
\eprop
We see that the short rate $r$ and the  Girsanov kernel $\vp$ has exactly the same structural 
form as the standard case formulas \cref{me10}-\cref{me11}. We now move to the  stochastic discount factor and after some calculations we have the following version of Theorem \ref{m1}.

\bprop \label{ned15}
The stochastic discount factor $M$ is determined by
\beq \label{ned16}
d \ln \pa{M_t}=d\ln g_x(t,X_t)
\eeq
where $g_x$ is evaluated at $(0,X_t)$.
Alternatively, we can write $M$ as
\beq \label{ned17}
M_t=U_c(\hat{c}_t)\cdot \exp{\krull{\int_0^t\frac{g_{xt}(0,X_s)}{g_x(0,X_s)}ds}}
\eeq
\eprop
We can also refer to Proposition \ref{mi12} and conclude that  the infinitesimal SDF is given by the formula
\beq \label{ned20}
\frac{M_{t+dt}}{M_t}=\frac{g_x\pa{dt,X_{t+dt}}}{g_x(0,X_t)},
\eeq
or, alternatively, by the formula
\beq \label{ned21}
\frac{M_{t+dt}}{M_t}=\frac{U_c\pa{\hat{c}_{t+dt}}}{U_c\pa{\hat{c}_t}}\cdot \frac{g_x\pa{dt,X_{t+dt}}}{g_x\pa{0,X_{t+dt}}}.
\eeq
\subsubsection{Log utility} \label{logut}
We now specialize to the case of log utility, so the utility functional has the form
\beq 
\EtqX{t,x}{}{\int_t^\infty \beta (\tau-t)\ln \pa{c_{\tau}}d\tau}.
\eeq
Given some experience from the standard time consistent case, we now make the Ansatz
\beq
g(t,x)=a_t \ln (x) + b_t,
\eeq
where $a$ and $b$ are deterministic functions of time. The natural boundary conditions are
\beq 
\lim_{t \rightarrow \infty} a_t=0,\quad \lim_{t \rightarrow \infty} b_t=0.
\eeq
With this Ansatz we have
\beq
g_t=\dot{a} \ln(x) + \dot{b}, \quad
g_x=\frac{a}{x},\quad
g_{xx}=- \frac{a}{x^2},
\eeq
so from the first order condition \cref{ned10}  for $c$ we  obtain
\beq
\hat{c}(x)=\frac{x}{a_0}.
\eeq
From Proposition \ref{ned5} we  obtain the short rate and the Girsanov kernel as
\beq
r=\alpha - \sigma^2, \quad \vp = - \sigma.
\eeq
The function $a$ and $b$ are determined by \cref{ned81}. We obtain
\beqarno
\dot{a} \ln (x) + \dot{b} +a \sigma^2 + \pa{\alpha - \sigma^2}a   -\frac{a}{a_0}- 
\frac{\sigma^2}{2} a+ \beta \ln (x) - \beta \ln \pa{a_0}=0
\eeqarno
We thus obtain the  ODE 
\beq
\dot{a}_t=-\beta(t),
\eeq
and with the boundary condition $\lim_{t \rightarrow \infty} a_t=0$ this gives us
\beq
a_t=\int_t^{\infty}\beta (s)ds.
\eeq
We also have an obvious ODE for $b$, but this is of little interest for us.

In order to determine the SDF we use Proposition \ref{ned15} and compute
$$
\frac{g_{xt}(0,x)}{g_x(0,x)}=\frac{\dot{a}_0}{a_0}=-\frac{1}{a_0}.
$$
and we have the following result.

\bprop
For the case of log utility, the  stochastic discount factor is given by
\beq
M_t=\frac{1}{X_t}\cdot e^{- t/a_0}.
\eeq
where 
$$
a_0=\int_0^{\infty}\beta (s)ds.
$$
\eprop

\subsubsection{Power utility} \label{po}
We now turn to the more complicate case of power utility so we have
$$
U(c)=\frac{c^{\gamma}}{\gamma}
$$ 
where $\gamma < 1$. We make the obvious Ansatz
\beq 
g(t,x)=a_t \frac{x^{\gamma}}{\gamma}
\eeq
We readily obtain
\beq
g_t=\dot{a}\frac{x^{\gamma}}{\gamma} ,\quad g_x=a  x^{\gamma -1},\quad g_{xx}=a(\gamma -1) x^{\gamma -2}, \quad g_{xt}=\dot{a} x^{\gamma -1}
\eeq
The first order condition for $c$ is
$$
c^{\gamma-1}=a_0 x^{\gamma -1}
$$
so the equilibrium consumption is given by
\beq
\hat{c}(x)=D x
\eeq
where
$$
D=a_0^{-1/(1-\gamma)}
$$
From Proposition \ref{ned5} we  obtain the short rate and the Girsanov kernel as
\beqarno
\vp&=&- \sigma (1- \gamma ),\\
r&=&\alpha - \sigma^2 (1-\gamma ).
\eeqarno
The function $a$ is again determined by \cref{ned81}. We obtain
\beqarno
\dot{a}\frac{x^{\gamma}}{\gamma}+ xa\sigma^2 (1-\gamma)x^{\gamma -1}+
xa\pa{r-D}x^{\gamma -1}+a \frac{\sigma^2}{2}x^2(\gamma -1) x^{\gamma -2} + \beta \frac{D^{\gamma}x^{\gamma}}{\gamma}=0
\eeqarno
which gives us the linear ODE 
\beqarno
\dot{a}_t+ Aa_t + B\beta (t)=0.
\eeqarno
with
\beqarno
A&=&\alpha -D - \frac{\sigma^2}{2}(1-\gamma),\\
B&=&\frac{D^{\gamma}}{\gamma}
\eeqarno
The function $a$ is thus given by
\beq
a_t=a_0e^{-At}-B\int_0^te^{-A(t-s)}\beta (s)ds.
\eeq
Using Proposition \ref{ned15} we thus have the following result.
\bprop
For the case of power utility the stochastic discount factor is given by
$$
M_t={X_t}^{\gamma - 1}e^{\frac{\dot{a}_0}{a_0}t}
$$ 
\eprop

From these examples with non exponential discounting we see that the risk free rate and Girsanov kernel only depend on the production opportunities in the economy. These objects are unaffected by the time inconsistency stemming from non-exponential discounting. Equilibrium consumption, however, is determined by the discounting function of the representative agent. In particular, we see that non exponential discounting has an effect on the marginal propensity to consume  from wealth and thus affects the equilibrium level of wealth in the economy.

\section{Conclusion and future research}
In this paper we have presented a fairly general class of time inconsistent stochastic control problems.
Using a game theoretic perspective we have derived a system of
equations for the determination of the subgame
perfect Nash equilibrium control, as well as for the corresponding equilibrium value function. The system is
an extension of the standard dynamic programming equation for time consistent problems. We have studied a couple of concrete examples, and in particular we have studied the effect of time inconsistency in the framework of general equilibrium for a production economy.

Some obvious open research problems are the following.
\bei
\item
In Section \ref{cpd} we informally derived the continuous time extended HJB system as a limit using a discretization argument.
It would obviously be valuable to have a theorem which rigorously proves convergence of the discrete time theory to the continuous time limit.
For the quadratic case, this is done in \cite{Czi}, but the general problem is completely open (and probably very hard).
\item
A related (hard) open problem is to prove existence and/or uniqueness for solutions of the extended HJB system.
\item
The present theory depends critically on the Markovian structure of the model. It would be interesting to see what can be
one without this assumption.
\item
The equilibrium model in Section \ref{eqp} can  be extended to a multidimensional model with several underlying factors. 
This is the subject of a forthcoming paper.
\ei



\end{document}